\input ./style/arxiv-general.cfg
\documentclass[aap,MSNbibl,seceqn,dvips]{arximspdf}
\makeatletter
   \@ifpackageloaded{graphicx}{}{\usepackage{graphicx}}
\makeatother

\doi{10.1214/14-AAP1068} 
\volume{25}
\issue{6}
\pubyear{2015}
\firstpage{3095}
\lastpage{3138}
\docsubty{FLA}

\makeatletter
\def\triangle{\Delta}
\def\mid{|}

\newcommand{\rrvert}{\vert}

\newcommand{\rrVert}{\Vert}
\newcommand{\llvert}{\vert}
\newcommand{\llVert}{\Vert}
\newcommand{\eqref}[1]{(\ref{#1})}
\newtheorem{theorem}{Theorem}[section]
\newtheorem{proposition}[theorem]{Proposition}
\newtheorem{lemma}[theorem]{Lemma}
\newtheorem{corollary}[theorem]{Corollary}
\newproclaim{remark}[theorem]{Remark}
\newproclaim{example}[theorem]{Example}
\newtheorem{guess}[theorem]{Hypothesis}
\makeatother

\begin{document}
\begin{frontmatter}

\title{A probabilistic interpretation of the parametrix method\thanksref{T1}}
\runtitle{A probabilistic interpretation of the parametrix method}
\thankstext{T1}{Supported by grants from the Japanese government.}

\begin{aug}
\author[A]{\fnms{Vlad}~\snm{Bally}\corref{}\ead[label=e1]{bally@univ-mlv.fr}}
\and
\author[B]{\fnms{Arturo}~\snm{Kohatsu-Higa}\ead
[label=e2]{arturokohatsu@gmail.com}}
\runauthor{V. Bally and A. Kohatsu-Higa}
\affiliation{Universit\'e Paris-Est Marne-la-Vall\'ee and INRIA, and
Ritsumeikan University and Japan Science and Technology Agency}
\address[A]{LAMA (UMR CNRS, UPEMLV, UPEC)\\
Universit\'e Paris-Est Marne-la-Vall\'ee\\
and\\
MathRisk\\
INRIA\\
F-77454 Marne-la-Vall\'ee\\
France\\
\printead{e1}}
\address[B]{Department of Mathematical Sciences\\
Ritsumeikan University\\
\quad and Japan Science and Technology Agency\\
1-1-1 Nojihigashi\\
Kusatsu, Shiga, 525-8577\\
Japan\\
\printead{e2}}
\end{aug}

\received{\smonth{6} \syear{2013}}
\revised{\smonth{6} \syear{2014}}

%
\begin{abstract}
In this article, we introduce the parametrix technique in order to
construct fundamental solutions as a general method based on semigroups
and their generators. This leads to a probabilistic interpretation of the
parametrix method that
is amenable to Monte Carlo simulation. We consider the explicit
examples of continuous diffusions and jump driven stochastic
differential equations with H\"older continuous coefficients.
\end{abstract}

%
\begin{keyword}[class=AMS]
\kwd[Primary ]{35K10}
\kwd{35K15}
\kwd{65C20}
\kwd[; secondary ]{65C05}
\kwd{65C30}
\end{keyword}
\begin{keyword}
\kwd{Parametrix}
\kwd{stochastic differential equations}
\kwd{density}
\kwd{Monte Carlo methods}
\end{keyword}
\end{frontmatter}

\section{Introduction}\label{sec1}

The parametrix technique for solving parabolic partial differential
equations (PDEs)
is a classical method in order to expand the fundamental solution of
such an
equation in terms of a basic function known as the parametrix. This is the
parallel of the Taylor expansion of a smooth function in terms of
polynomials.

The concept of order of the polynomial in the classical Taylor
expansion is replaced by multiple integrals whose order increases as
the expansion becomes more accurate. This method has been successfully
applied to many
equations and various situations. Its success is due to its flexibility
as it can be applied to a wide variety of
PDEs. It has been successfully extended to other
situations for theoretical goals (see, e.g., \cite{KM,KM1,S3,S2} and \cite{S1}). In \cite{F}, the authors consider
the parametrix as an analytical method for approximations for
continuous diffusions.
These analytical approximations may be used as deterministic
approximations and are highly accurate in the cases where the sum
converges rapidly. In general, higher order integrals are difficult to
compute and, therefore, this becomes a limitation of the method.

The goal of the present paper is to introduce a general probabilistic
interpretation
of the parametrix method based on semigroups, which not only
reexpresses the arguments of the
method in probabilistic terms, but also to introduce an alternative
method of
simulation with no approximation error.
This leads to the natural emergence
of the difference between the generators of the process and its
parametrix in the same manner as the concept of derivative appears in
the classical Taylor expansion.

Let us explain the above statement in detail.
The first step in the Monte Carlo approach for approximating the
solution of the parabolic partial differential equation $\partial
_{t}u=Lu$ is to construct the Euler scheme which
approximates the continuous diffusion process with infinitesimal
operator $L$. To fix the
ideas, consider the diffusion process $X_{t}\equiv X_t(x)$ solution of
the following stochastic differential equation (SDE):
%
\begin{equation}
dX_{t}=\sum_{j=1}^{m}
\sigma_{j}(X_{t})\,dW_{t}^{j}+b(X_{t})\,dt,\qquad
t\in [0,T], X_0=x_0, \label{eq:1.1}
\end{equation}
where $W$ is a multidimensional Brownian motion and $\sigma
_{j},b\dvtx \mathbb{R}^{d}\rightarrow\mathbb{R}^{d}$
are smooth functions. We denote by $P_{t}f(x)=\mathbb{E}[f(X_{t}(x))]$ the
semigroup associated to this diffusion process. The infinitesimal
generator associated to
$P$ is defined for $f\in C^2_c(\mathbb{R}^d)$ as
%
\begin{equation}
Lf(x)=\frac{1}2 \sum_{i,j}a^{i,j}(x)
\partial _{i,j}^{2}f(x)+b^{i}(x)\partial
_{i}f(x), \qquad a:=\sigma\sigma^{\ast}. \label{Int2}
\end{equation}
By the Feynman--Kac formula, one knows that $u(t,x):=P_tf(x)$ is the
unique solution to $\partial_{t}u=Lu$ satisfying the initial condition
$u(0,x)=f(x)$.
Therefore, the goal is to approximate $X$ first and then the
expectation in $P_tf(x)={\mathbb{E}}[f(X_{t}(x))]$
using the law of large numbers which leads to the Monte Carlo method.

Now, we describe some stochastic approximation methods for $X$.
Given a partition of $[0,T]$, $\pi=\{0=t_0<\cdots<t_n=T\} $, the Euler
scheme associated to this time grid is defined as $X^\pi_0(x)=x$
%
\begin{eqnarray}\label{Int3}
X_{t_{k+1}}^{\pi}(x)&=&X_{t_k}^{\pi}(x)+\sum
_{j=1}^{m}\sigma _{j}
\bigl(X_{t_k}^{\pi}(x)\bigr) \bigl(W^j_{t_{k+1}}-W^j_{t_{k}}
\bigr)
\nonumber
\\[-8pt]
\\[-8pt]
\nonumber
&&{}+b\bigl(X_{t_{k}}^{\pi
}(x)\bigr) (t_{k+1}-t_k).
\end{eqnarray}
It is well known (see \cite{TT}) that $X^\pi\equiv X^\pi(x)$ is an
approximation scheme of $X$ of order one. That is, there exists a
constant $C_f(x)$ such that
%
\begin{eqnarray}\label{Int4}
&&\bigl\llvert \mathbb{E}\bigl[f\bigl(X_{T}(x)\bigr)\bigr]-\mathbb{E}
\bigl[f\bigl(X_{T}^{\pi}(x)\bigr)\bigr]\bigr\rrvert
\nonumber
\\[-8pt]
\\[-8pt]
\nonumber
&&\qquad\leq
C_{f}(x)\max\{t_{i+1}-t_i;i=0,\ldots,n-1\}
\end{eqnarray}
for $f$ measurable and bounded (see \cite{bally96a,bally96b})
and under
strong regularity assumptions on the coefficients $\sigma_{j}$ and $b$.


Roughly
speaking, the parametrix method is a deterministic method with the following
intuitive background: in short time the diffusion $X_{t}(x_{0})$ is
close to
the diffusion with coefficients ``frozen'' in the starting point
$x_{0}$. So one
may replace the operator $L$ by the operator $L^{x_{0}}$ defined as
\[
L^{x_{0}}f(x)=\frac{1}2 \sum_{i,j}a^{i,j}(x_{0})
\partial _{i,j}^{2}f(x)+b^{i}(x_{0})
\partial_{i}f(x)
\]
and one may replace the semigroup $P_{t}$ by the semigoup $%
P_{t}^{x_{0}}$ associated to $L^{x_{0}}$. Clearly, this is
the same idea as the one which leads to the construction of the Euler scheme
(\ref{Int3}). In fact, notice that the generator of the one step (i.e.,
$\pi=\{0,T\}$) Euler scheme $X^\pi(x_0)$ is given by $L^{x_0}$.

The goal of the present article is to give a probabilistic
representation formula based on the parametrix method. This formula
will lead to simulation procedures for $\mathbb{E}[f(X_T)]$ with no
approximation error which are based on the weighted sample average of
Euler schemes with random partition points given by the jump times of
an independent Poisson process.
In fact, the first probabilistic representation formula (forward
formula) we intend to prove is the following:
%
\begin{equation}
\label{eq:ESF} \mathbb{E}\bigl[f(X_{T})\bigr]=e^{T}
\mathbb{E} \Biggl[f\bigl(X_{T}^{\pi}\bigr)\prod
_{k=0}^{J_{T}-1}\theta_{\tau_{k+1}-\tau_k}
\bigl(X_{\tau
_{k}}^{\pi},X_{\tau_{k+1}}^{\pi}\bigr)
\Biggr].
\end{equation}
Here, $\tau_0=0$ and $\pi:=\{\tau_{k},k\in\mathbb{N}\}$, are the jump
times of a
Poisson process $\{J_{t};t\in[0,T]\}$ of parameter one and $%
X^{\pi}$ is a continuous time Markov process which satisfies
\[
P\bigl(X_{\tau_{k+1}}^{\pi}\in dy\mid\{\tau_k,k\in
\mathbb{N}\}, X_{\tau
_{k}}^{\pi}=x_{0}
\bigr)=P_{\tau_{k+1}-\tau
_{k}}^{x_{0}}(x_{0},dy)=p^{x_0}_{\tau_{k+1}-\tau
_{k}}(x_0,y)\,dy.
\]
In the particular case discussed above, then $X^\pi$ corresponds in
fact to an Euler scheme with random partition points.
$\theta_t\dvtx \mathbb{R}^d\times\mathbb{R}^d\rightarrow\mathbb{R}$ is a
weight function to be described later.

Before discussing the nature of the above probabilistic representation
formula, let us remark that a similar formula is available for
one-dimensional diffusion processes (see \cite{BR}) which is strongly
based on explicit formulas that one may obtain using the Lamperti
formula. Although many elements may be common between these two
formulations, the one presented here is different in nature.

In order to motivate the above formula \eqref{eq:ESF}, let us give the
following basic heuristic argument that leads to the forward parametrix method:
%
\begin{eqnarray}
\label{eq:2.1} P_tf(x)-P_t^{x}f(x)&=&\int
_0^t\partial_s
\bigl(P^x_{t-s}P_sf\bigr) (x)\,dt
\nonumber
\\[-8pt]
\\[-8pt]
\nonumber
&=&\int
_0^tP^x_{t-s}
\bigl(L-L^{x}\bigr)P_sf(x)\,ds.
\end{eqnarray}
Here, we suppose that $P_sf\in \operatorname{Dom}(L-L^x)$. The above
expression is already an equivalent of a Taylor expansion of order one
where the notion of first-order derivative is being replaced by
$L-L^{x}$. Its iteration will lead to higher order Taylor expansions.
Another way of looking at this is to consider \eqref{eq:2.1} as a
Volterra equation in $Pf(x)$. This will be our point of view here.

In fact, if $t$ is considered as the time variable and considering
\eqref{eq:2.1} as an equation, one sees, after an application of the
integration by parts on the diffferential operator $L-L^x$, that $\{
P_tf; t\in[0,T]\}$ for $f\in C^\infty_c$ can also be considered as a
solution of the following Volterra type linear equation:
%
\begin{eqnarray}
\label{eq:3.1} P_tf(x)&=&P_t^xf(x)+\int
_0^t\int p_{t-s}^x(x,y_1)
\bigl(L-L^x\bigr)P_sf(y_1)\,dy_1\,ds
\nonumber
\\[-8pt]
\\[-8pt]
\nonumber
&=&P_t^xf(x)+\int_0^t
\int\theta _{t-s}(x,y_1)p_{t-s}^x(x,y_1)P_sf(y_1)\,dy_1\,ds.
\end{eqnarray}
Here, we have used the function $\theta$, defined as
%
\begin{eqnarray}\label{eq:3.11}
&&\theta_{t-s} (x,y_1)p^z_{t-s}(x,y_1)
|_{z=x} \nonumber\\
&&\qquad
=\bigl(L-L^z\bigr)^*p_{t-s}^z(x,
\cdot) (y_1) |_{z=x}
\nonumber
\\[-8pt]
\\[-8pt]
\nonumber
&&\qquad=\frac{1}2\sum_{i,j}
\partial_{i,j}\bigl( \bigl(a^{i,j}(\cdot)-a^{i,j}(z)
\bigr)p^z_{t-s}(x,\cdot)\bigr) (y_1)\\
&&\qquad\quad{}-\sum
_i\partial _i\bigl(
\bigl(b^i(\cdot)-b^i(z)\bigr)p^z_{t-s}(x,
\cdot)\bigr) (y_1) \bigg|_{z=x}.\nonumber
\end{eqnarray}
Equation \eqref{eq:3.1} can be iterated in $P_tf$ in order to obtain a
series expansion of the solution which is
the equivalent of the Taylor expansion of $P_tf$.

We may note that the second term in this expansion for $t=T$ can be
rewritten using the Euler scheme with partition points $\pi=\{0,T-s,T\}
$ as
%
\begin{eqnarray}
\label{eq:2.11} &&\int_0^T\int
\theta_{T-s}(x,y_1)p_{T-s}^x(x,\cdot)
(y_1)P^{y_1}_sf(y_1)\,dy_1\,ds
\nonumber
\\[-8pt]
\\[-8pt]
\nonumber
&&\qquad
=\int_0^T\mathbb{E}\bigl[f
\bigl(X_T^{\pi}(x)\bigr)\theta_{T-s}
\bigl(x,X_{T-s}^{\pi
}(x)\bigr)\bigr]\,ds.
\end{eqnarray}
This is the first step toward the construction of what we call the
forward parametrix method.\setcounter{footnote}{1}\footnote{This also explains the logic behind
the choice of variables in the integrals. Through the rest of the paper
$y_i$ will denote the integrating variables in the order given by the
corresponding Euler scheme. Similar rule will apply with the times
$t_i$, $i\in\mathbb{N}$. For example, we will have that $t_1=t-s_1$ in
\eqref{eq:2.1}.} It requires the regularity of the coefficients and it
is based on the usual Euler scheme for the sde
\eqref{eq:1.1}. Note that \eqref{eq:2.11} will be associated with the
term $J_T=1$ in \eqref{eq:ESF}. In fact, if there is only one jump of
the Poisson process $J$ in the interval $[0,T]$, then the distribution
of the jump is uniform in the interval $[0,T]$. This leads to the
probabilistic interpretation of the time integral in \eqref{eq:2.11}.

Let us now discuss an alternative to the above method which requires
less regularity conditions
on the coefficients of \eqref{eq:1.1}. This method will be called the
backward parametrix method and it is obtained by duality arguments as
follows. That is, consider for two functions $f,g\in C^\infty_c(\mathbb
{R}^d)$ the pairing $\langle f, P_t^*g \rangle$. Then we will use the
approximating semigroup $\hat{Q}_tg(y):=(P_t^y)^*g(y)=\int
g(x)p_{t}^{y}(x,y)\,dx$. A~similar heuristic argument gives
%
\begin{eqnarray}
\label{eq:3.2} P_t^*g(y)-\bigl(P_t^y
\bigr)^*g(y)&=&\int_0^t\int
\bigl(L-L^y\bigr)p_{t-s}^y(\cdot ,y)
(y_1)P^*_sg(y_1)\,dy_1\,ds
\nonumber
\\[-8pt]
\\[-8pt]
\nonumber
&=&\int_0^t\int\hat\theta
_{t-s}(y_1,y)p_{t-s}^y(y_1,y)P^*_sg(y_1)\,dy_1\,ds.
\end{eqnarray}
Note that in this case the operator $L-L^y$ is applied to the density
function $p_{t-s}^y(\cdot,y)$ with the coefficients frozen at $y$,
therefore, no derivative of the coefficients is needed in this approach.
In fact,
%
\begin{eqnarray}
\label{eq:3.21} &&\hat\theta_{t-s}(y_1,y)p_{t-s}^z(y_1,y)
|_{z=y}
\nonumber
\\[-8pt]
\\[-8pt]
\nonumber
&&\qquad= \bigl(L-L^z\bigr)p_{t-s}^z(
\cdot,y) (y_1) |_{z=y}
\\
&&\qquad= \biggl\{\frac{1}2 \sum_{i,j}
\bigl(a^{i,j}(y_1)-a^{i,j}(z)\bigr)
\partial_{i,j}
\nonumber
\\[-8pt]
\\[-8pt]
\nonumber
&&\hspace*{32pt}\qquad{}+\sum_i\bigl(b^i(y_1)-b^i(z)
\bigr)\partial _i \biggr\} p^z_t(\cdot,y)
(y_1) \bigg|_{z=y}.
\end{eqnarray}
As before, we can obtain a probabilistic representation. In this case,
one has to be careful with the time direction. In fact, due to the
symmetry of the density function $p_{t-s}^y(y_1,y)$ one interprets it
as the density of the Euler scheme at $y_1$ started at $y$. Therefore,
the sign of the drift has to be changed leading to what we call the
backward parametrix method. In the particular case that $f$ is a
density function, it will be interpreted as a ``backward running'' Euler
scheme from $T$ to $0$ with random initial point with density $f$. The
test function $g$ is replaced by a Dirac delta at the initial point of
the diffusion $x_0$.
See Section~\ref{sec:back} for precise statements.


Therefore, the behavior of forward and backward methods are different.
In fact, the forward method applies when the coefficients are regular.
In many applied situations, one may have coefficients which are just H\"
older continuous and, therefore, the forward method does not apply. In
that case, one may apply the backward method which demands less regularity.
For this reason, the treatment of the forward method and the backward
method are essentially different and they are treated separately.
Issues related to simulation will be discussed in another article.

Our article is structured as follows: In Section~\ref{sec:2}, we give
the notation used throughout the paper. In Section~\ref{sec:FE}, we
discuss the existence, uniqueness and regularity properties of the
solution of the linear Volterra equations of the type \eqref{eq:3.1} or
\eqref{eq:3.2} which will be applicable to both probabilistic
representation formulas to be discussed later. In Section~\ref{sec:3},
we provide a general abstract framework based on semigroups for which
our two approaches (forward and backward) can be applied. The main
hypotheses applying to both methods are given in this section. In
Section~\ref{sec:forward}, we give the analytical form of the forward
method. In Section~\ref{sec:probaforw}, we give the probabilistic
representation, and in Section~\ref{sec:rgul}, we give the
continuity and differentiability properties of the density functions.
This is the usual application of the parametrix method.

In Section~\ref{sec:back}, we give the backward approach which was
first introduced in \cite{MS}.
We also give the probabilistic interpretation and the regularity
results corresponding to the backward method in parallel sections.

In Section~\ref{sec:contdiff}, we consider our main examples. The first
corresponds to the continuous diffusion with uniformly elliptic
diffusion coefficient. We see in Section~\ref{sec:6.1} that in the
forward approach we need the coefficients to be smooth. While in
Section~\ref{sec:6.2}, we show that in order for the backward
approach to be applicable, we only require the coefficients to be H\"
older continuous. In Section~\ref{sec:jumps}, we also consider the
case of a jump driven SDE where the L\'evy measure is of stable type in
a neighborhood of zero. This example is given with two probabilistic
interpretations.

We close with some conclusions: an \hyperref[app]{Appendix} and the References section.

\section{Some notation and general definitions}
\label{sec:2}
We now give some basic notation and definitions used through this article.
For a sequence of operators $S_{i}$, $i=1,\ldots,n$ which do not necessarily
commute, we define $\prod_{i=1}^{n}S_{i}=S_{1}\cdots S_{n}$ and $\prod_{i=n}^1S_{i}=S_n\cdots S_1$. We will denote by $I$, the
identity matrix or identity operator and $S^{\ast}$ will denote the adjoint
operator of $S$. $\operatorname{Dom}(S)$ denotes the domain of the operator $S$. If
the operator $S$ is of integral type, we
will denote its associated measure $S(x,dy)$ so that $Sf(x)=\int
f(y)S(x,dy)$. All space integrals will be taken over $\mathbb{R}^d$.
For this reason, we do not write the region of integration which we
suppose clearly understood.
Also in order to avoid long statements, we may refrain from writing
often where the time and space variables take values supposing that
they are well understood from the context.

In
general, indexed products where the upper limit is negative are
defined as 1 or $I$. In a similar fashion, indexed sums where the upper
limit is negative are defined as zero.

As it is usual, $A\leq B$ for two matrices $A$ and $B$, denote the fact that
$A-B$ is positive definite. Components of vectors or matrices are
denoted by
superscript letters. When the context makes it clear we denote by
$\partial
_{i}f$ the partial derivative operator with respect to the $i$th
variable of the function $f$ and similarly for higher order
derivatives. For example,
derivatives with respect to a multi-index $\beta$, of length $\llvert
\beta
\rrvert  $, are denoted by $\partial_{\beta}f$. Time derivatives will
be denoted by $\partial_t$.

We denote by $\delta_{a}(dx)$ the point mass measure concentrated in
$\{a\}$, $B(x,r)$ denotes the ball of center $x\in\mathbb{R}^d$ and
radius $r>0$, $%
[x]$ denotes the ceiling or smallest integer function for $x\in\mathbb
{R}$ and $\mathbb{R}_+\equiv
(0,\infty)$. The indicator function of the set $A$ is denoted by
$1_A(x)$, $C(A)$ denotes the space of real valued functions continuous
in the set $A$. The space of real valued measurable bounded functions
defined on $A$ is denoted by $L^\infty(A)$. Similarly, the space of
continuous bounded functions
in $A$ is denoted by $C_{b}(A)$. The space of real valued infinitely
differentiable functions with compact support defined on $\mathbb{R}^d$
is denoted by $C_{c}^{\infty}(\mathbb{R}^{d})$. The space of $\mathbb
{R}^{l}$-valued bounded
functions defined on $\mathbb{R}^d$ with bounded derivatives up to
order $k$ is denoted by $%
C_{b}^{k}(\mathbb{R}^d;\mathbb{R}^{l})$. The
norm in this space is defined as $\llVert  f\rrVert  _{k,\infty
}=\max_{\llvert  \beta\rrvert  \leq k}\sup_{x\in\mathbb{R}^d}\llvert
\partial_{\beta
}f(x)\rrvert  $. In the particular case that $k=0$, we also use the
simplified notation
$\|f\|_\infty\equiv\|f\|_{0,\infty}$.

The
multidimensional Gaussian density at $y\in\mathbb{R}^d$ with mean zero and
covariance matrix given by the positive definite matrix $a$ is denoted
by $q_a(y)$.
Sometimes we abuse the notation denoting by $q_{t}(y)$, for $y\in
\mathbb{R}^d$, $%
t>0$ the Gaussian density corresponding to the variance--covariance
matrix $%
tI $. Similarly, $H^i_{a}(y)$ and $H^{i,j}_a(y)$ for $i, j\in\{1,\ldots,d\}
$, denote the multidimensional version of the
Hermite polynomials of order one and two. Exact definitions and some of
the properties of Gaussian densities used
throughout the article are given in Section~\ref{sec:Gaussian}.

Constants will be denoted by $C$ or $c$, we will not give the explicit
dependence on parameters of the problem unless it is needed in the
discussion. As it is usual, constants may change from one line to the next
although the same symbol may be used.

In the notation throughout the article, we try to denote by $x$ the
starting point of the diffusion and $y$ the arrival point with $z$
being the parameter value where the operator $L^z$ is frozen at. In the
forward method, $z$ will be the starting point $x$ and in the backward
method $z$ will be the arrival point $y$.
Due to the iteration procedure, many intermediate points will appear
which will be generally denoted by $y_i$, $i=0,\ldots,n$, always going
from $y_0=x$ toward $y_{n}=y$ in the forward method and from $y_0=y$ to
$y_n=x$ in the backward method. As stated previously, the time
variables will be evolving forward in the sense of the Euler scheme
if they are denoted by~$t_i$, $i=0,\ldots,n$ from $t_0=0$ to $t_n=t$ or
backwardly if denoted by $s_i$, $i=0,\ldots,n$ from $s_0=t$ to $s_n=0$.

\section{A functional linear equation}
\label{sec:FE}
In this section, we consider a functional equation of Volterra type
which will include both equations \eqref{eq:3.1} and \eqref{eq:3.2}.
Therefore, this represents and abstract framework which includes the
forward and backward method.

We consider a jointly measurable functions $a \dvtx (0,T]\times\mathbb
{R}^d\times
\mathbb{R}^d\rightarrow\mathbb{R}$ and we define the operator
\[
U_{a}f(t,x)=\int_{0}^{t}\int
_{\mathbb{R}^{d}}f(s,y)a_{t-s}(x,y)\,dy\,ds.
\]
Our aim is to solve the equation
%
\begin{equation}
f=g+U_{a}f \label{eq:app1}
\end{equation}
and to study the regularity of the solution. Formally, the unique
solution is obtained by iteration and given
by $H_{a}g:=\sum_{n=0}^{\infty}U_{a}^{n}g$. In order to make this
calculation mathematically sound, we have to study the
convergence of the series.
For this, we consider the iterations of the operator $U_{a}$. We define
$U_{a}^{0}$ to be
the identity operator, $U_{a}^{1}=U_{a}$ and we define by recurrence $%
U_{a}^{n}=U_{a}^{n-1}U_{a}$.

\begin{lemma}
If $I_{a}(t,x):=\int_{\mathbb{R}^{d}}\llvert
a_{t}(x,y)\rrvert  \,dy\in L^\infty([0,T]\times\mathbb{R}^d)$ and
$g\in \break L^\infty([0,T]\times\mathbb{R}^d)$, then the equation \eqref
{eq:app1} has a unique solution in the space $\mathcal{A}:=\{f\in
L^\infty([0,T]\times\mathbb{R}^d);\lim_{N\rightarrow\infty}\|U_a^Nf\|
_\infty=0\}$.
\end{lemma}

\begin{pf} In fact,
\[
\bigl\llVert U_{a}g(t,\cdot)\bigr\rrVert _{\infty}\le\int
_0^t\bigl\Vert g(s,\cdot)\bigr\Vert_{\infty}\bigl
\llVert I_{a}(t-s,\cdot)\bigr\rrVert _{\infty}\,ds.
\]
Then by induction it follows that
\[
\bigl\llVert U_{a}^{n}g\bigr\rrVert _{\infty}\leq
\frac{T^{n}\llVert
I_{a}\rrVert _{\infty}^{n}}{n!}\llVert g\rrVert _{\infty}.
\]
This means that the infinite sum $H_{a}g:=\sum_{n=0}^{\infty
}U_{a}^{n}g$ converges absolutely and, therefore, is well defined in
$L^\infty([0,T]\times\mathbb{R}^d)$.
Furthermore, it is easy to see that the sum is a solution of equation
\eqref{eq:app1} satisfying
$H_ag\in\mathcal{A}$.

For any solution $f$ of \eqref{eq:app1}, one obtains by iteration that
\[
f=\sum_{n=0}^{N }U_{a}^{n}g+U^N_af.
\]
Therefore, if $f\in\mathcal{A}$ then $f$ satisfies $f=H_ag$. From here,
one obtains the uniqueness.
\end{pf}
Unfortunately, in our case $\int|a_{t}(x,y)|\,dy$ blows up as
$t\rightarrow0$ and we center our discussion on this matter. We see
from \eqref{eq:3.11} and \eqref{eq:3.21} that the rate of divergence is
determined by the regularity of the coefficients. We will call this
regularity index $\rho$ in what follows.
In order to introduce our main assumption, we define for a function
$\beta\dvtx (0,T]\times\mathbb{R}^d\times\mathbb{R}^d\rightarrow
\mathbb{R}_+$, the class of functions $\Gamma_\beta$ such that
there exists a positive constant $C$ which satisfies the following
inequality for every $n\in\mathbb{N},y_{0},y_{n+1}\in\mathbb{R}^d$
and every $%
\delta_{i}>0,i=1,\ldots,n$ with $s(\delta):=\sum_{i=1}^{n}\delta
_{i}$,
%
\begin{equation}
\label{eq:h2} \int dy_{1}\cdots\int dy_{n}\prod
_{i=0}^{n}\gamma _{\delta_{i}}(y_{i},y_{i+1})
\leq C^{n+1}\beta_{s(\delta)}(y_{0},y_{n+1}).
\end{equation}

\begin{guess}\label{(h1)}
%
There exists a positive constant $C$ and a function
$%
\gamma\in\Gamma_\beta$ such that $\sup_{(t,x)\in[0,T]\times\mathbb
{R}^d}\int|\gamma_t(x,y)|\,dy<\infty$ and $\rho
\in[0,1)$ such that for every $(t,x,y)\in(0,T]\times\mathbb
{R}^d\times\mathbb{R}^d$
%
\begin{equation}
\label{eq:h1} \bigl\llvert a_{t}(x,y)\bigr\rrvert \leq
\frac{C}{t^{\rho}}\gamma _{t}(x,y).
\end{equation}
Furthermore, there exists a function $\beta$ such that $\gamma\in\Gamma
_\beta$.
%
\end{guess}

If we define $
C_{\beta}(t,x):=\int_{0}^{t}\,ds\int dy\,\beta
_{t-s}(x,y)$ for $(t,x)\in[0,T]\times\mathbb{R}^d$, we will also
assume that

\begin{guess}
\label{(h3)}
%
$C_\beta$ is a bounded function.
%
\end{guess}

We denote\footnote{Notice that according to our remark about the
meaning of the variables in the integrals of \eqref{eq:2.11}, this
order of time is reversed with respect to the order in space.} $D_{T}=\{
((t,x),(s,y))\dvtx 0<s\leq t\leq T, x,y\in\mathbb{R}^d\}$. To the function~$a$, we
associate the function $A\dvtx D_{T}\rightarrow\mathbb{R}$
defined by
\[
A\bigl((t,x),(s,y)\bigr)=a_{t-s}(x,y).
\]
Then we define the operator $U_{a}\dvtx L^\infty([0,T]\times
\mathbb{R}^d)\rightarrow L^\infty([0,T]\times\mathbb{R}^d)$ by
\[
U_{a}f(t,x)=\int_{0}^{t}\,ds\int
dy\,f(s,y)A\bigl((t,x),(s,y)\bigr)=\int_{0}^{t}\,ds%
\int dy\,f(s,y)a_{t-s}(x,y).
\]
In fact, note that if $f\in L^\infty([0,T]\times
\mathbb{R}^d)$ then
\[
\bigl|U_{a}f(t,x)\bigr|\le\|f\|_\infty \int_{0}^{t}\,ds%
\int dy\frac{|\gamma_{t-s}(x,y)|}{(t-s)^\rho}\le C\|f\|_\infty t^{1-\rho}.
\]
Note in particular that this estimate implies that $U_a^nf$ is well
defined for $f\in L^\infty([0,T]\times
\mathbb{R}^d)$.
We also define for $0<s\leq t\leq T$ and $x,y\in
\mathbb{R}^d $
%
\begin{eqnarray}\label{J1}
A_{1}\bigl((t,x),(s,y)\bigr) &=&A\bigl((t,x),(s,y)\bigr),
\nonumber
\\
\qquad A_{n}\bigl((t,x),(s,y)\bigr) &=&\int_{s}^{s_{0}}\,ds_{1}
\int dy_{1}\int_{s_1}^{s_0}\,ds_{2}%
\int dy_{2}\cdots\int_{s_{n-2}}^{s_0}\,ds_{n-1}
\int dy_{n-1}\\
&&{}\times\prod_{i=0}^{n-1}A
\bigl((s_{i},y_{i}),(s_{i+1},y_{i+1})
\bigr)\nonumber
\end{eqnarray}
with the convention that $s_{0}=t,y_{0}=x$ and $s_{n}=s,y_{n}=y$, $n\ge
2$. Notice
that we have that $A_n$ is finite. That is, we have the following.

\begin{lemma}
\label{lemma:1.5}Assume Hypothesis \ref{(h1)}, then there exists a
constant $C(T,\rho)$ such that
%
\begin{eqnarray}
\label{J4} \bigl\llvert A_{n}\bigl((t,x),(s,y)\bigr)\bigr\rrvert
&\le& C \beta_{t-s}(x,y)\times\frac{C^{n}(T,\rho)}{[1+n\rho]!},
\\
\qquad A_{n+1}\bigl((t,x),(s,y)\bigr)&=&\int_{s}^{t}\,ds_1%
\int dy_1\,A_{n}\bigl((t,x),(s_1,y_1)
\bigr)A\bigl((s_1,y_1),(s,y)\bigr). \label{J2}
\end{eqnarray}
\end{lemma}

\begin{pf}
We use Hypothesis \ref{(h1)} and we obtain (with $%
s_{0}=t,y_{0}=x,s_{n}=s,y_{n}=y)$
\begin{eqnarray*}
&&\bigl\llvert A_{n}\bigl((t,x),(s,y)\bigr)\bigr\rrvert\\
&&\qquad\leq \int
_{s}^{s_{0}}\,ds_{1}\int dy_{1}
\int_{s_{1}}^{s_0}\,ds_{2}%
\int
dy_{1}\cdots\int_{s_{n-2}}^{s_0}\,ds_{n-1}
\int dy_{n-1}%
\\
&&\qquad\quad{}\times\prod_{i=0}^{n-1}(s_{i}-s_{i+1})^{-\rho}
\gamma_{s_{i}-s_{i+1}}(y_{i},y_{i+1})
\\
&&\qquad\leq C^n\beta _{s_{0}-s_{n}}(y_{0},y_{n})
\int_{s}^{s_{0}}\,ds_{1}\int
_{s_{1}}^{s_0}\,ds_{2}\cdots%
\int
_{s_{n-2}}^{s_0}\,ds_{n-1}\prod
_{i=0}^{n-1}(s_{i}-s_{i+1})^{-\rho}
\\
&&\qquad\le C^n\beta_{t-s}(x,y) (t-s)^{n(1-\rho) }
\frac{\Gamma^{n}(\rho
)}{[1+n\rho
]!}=\beta_{t-s}(x,y) \frac{C^{n}(T,\rho)}{[1+n\rho]!}
\end{eqnarray*}
the last inequality being a consequence of the change of variables
$s_i=s_0-t_i$ and Lemma~\ref{lem:app1} where we have set $C(T,\rho
)=CT^{1-\rho}\Gamma(\rho)$.
\end{pf}
Now that $A_n$ is well defined we can now give an
explicit formula for $U^n_a$.

\begin{lemma}Assume Hypotheses \ref{(h1)} and \ref{(h3)}. Let $f\in
L^\infty([0,T]\times\mathbb{R}^d)$ then
%
\begin{equation}
U_{a}^{n}f(t,x)=\int_{0}^{t}\,ds
\int dyA_{n}\bigl((t,x),(s,y)\bigr)f(s,y). \label{J3}
\end{equation}
\end{lemma}

\begin{pf}
For $n=1$, this is true by the definition of $U_{a}$. Suppose that this
is true
for $n$ and let us prove it for $n+1$. By (\ref{J2}),
\begin{eqnarray*}
&&U_{a}^{n+1}f(t,x)\\
&&\qquad=U_{a}^{n}U_{a}f(t,x)=
\int_{0}^{t}\,du%
\int dzA_{n}
\bigl((t,x),(u,z)\bigr)U_{a}f(u,z)
\\
&&\qquad=\int_{0}^{t}\,du\int dzA_{n}
\bigl((t,x),(u,z)\bigr)\int_{0}^{u}\,dv%
\int dwA\bigl((u,z),(v,w)\bigr)f(v,w)
\\
&&\qquad=\int_{0}^{t}\,dv\int dwf(v,w)\int
_{v}^{t}\,du%
\int dzA_{n}
\bigl((t,x),(u,z)\bigr)A\bigl((u,z),(v,w)\bigr)
\\
&&\qquad=\int_{0}^{t}\,dv\int dwf(v,w)\int
dzA_{n+1}\bigl((t,x),(v,w)\bigr).
\end{eqnarray*}
So (\ref{J3}) is proved. The integrability of the above expressions
follows from Lemma~\ref{lemma:1.5} and Hypothesis \ref{(h3)}.
\end{pf}

The main estimate in this section is the following. For this, we define
\[
C_{T}(\rho )=\sum_{n=0}^{\infty}
\frac{C^{n}(T,\rho)}{[1+n\rho]!}.
\]

\begin{theorem}
\label{lem:1.7}
\textup{(A)} Assume that Hypotheses \ref{(h1)} and \ref{(h3)} hold true.
Then the series
%
\begin{equation}
\label{eq:71} S_{a}\bigl((t,x),(s,y)\bigr)=\sum
_{n=1}^{\infty}A_{n}\bigl((t,x),(s,y)\bigr)
\end{equation}
is absolutely convergent and
\[
\bigl\llvert S_{a}\bigl((t,x),(s,y)\bigr)\bigr\rrvert \leq
C_{T}(\rho)\beta_{t-s}(x,y).
\]
\textup{(B)} Moreover, for $f\in L^\infty([0,T]\times\mathbb{R}^d)$ the
series
\[
H_{a}f(t,x):=\sum_{n=0}^{\infty}U_{a}^{n}f(t,x)
\]
is absolutely convergent and
\[
\bigl\llvert H_{a}f(t,x)\bigr\rrvert \leq C_{T}(
\rho)C_{\beta}(t,x)\llVert f\rrVert _{\infty}.
\]
Finally,
\[
H_{a}f(t,x)=\int_{0}^{t}\int
f(s,y)S_{a}\bigl((t,x),(s,y)\bigr)\,dy\,ds.
\]
\textup{(C)} Let $f,g\in L^\infty([0,T]\times\mathbb{R}^d)$ such that
\[
f=g+U_{a}f.
\]
Then $f=H_{a}g$.
\end{theorem}

\begin{pf}
From Lemma~\ref{lemma:1.5}, \eqref{J4}, we conclude that the
series $S_{a}((t,x),\break  (s,y))=\sum_{n=1}^{\infty}A_{n}((t,x),(s,y))$ is
absolutely convergent and $\llvert  S_{a}((t,x),\break (s,y))\rrvert \leq 
C_{T}(\rho)\beta_{t-s}(x,y)$.
We consider now a function $f\in L^\infty([0,T]\times\mathbb{R}^d)$.
As a
consequence of (\ref{J4}),
\[
U_{a}^{n}f(t,x)\leq\frac{C^{n}(T,\rho)}{[1+n\rho]!}\llVert f\rrVert
_{\infty}\int_{0}^{t}\int dy\,
\beta_{t-s}(x,y)\leq\frac{%
C^{n}(T,\rho)}{[1+n\rho]!}C_{\beta}(t,x)\llVert f
\rrVert _{\infty
}.
\]
It follows that the series $H_{a}f(t,x):=\sum_{n=0}^{\infty}U_{a}^{n}f(t,x)$
is absolutely convergent and $\llvert  H_{a}f(t,x)\rrvert \leq
C_{T}(\rho)C_{\beta}(t,x)\llVert  f\rrVert _{\infty}$.
Furthermore, from the above estimates it is clear that we can exchange
integrals and sums in order to prove that $H_ag$ is a solution to the
equation \eqref{eq:app1}. For any given bounded solution $f$ to \eqref
{eq:app1}, we obtain by iteration of the equation that the solution has
to be $H_ag$ and, therefore, we get the uniqueness.
\end{pf}
%
We give now a corollary with the study of the fundamental solution.
This will be used in order to obtain the density functions
corresponding to the operators appearing in \eqref{eq:3.1} and \eqref
{eq:3.2}. The proof follows directly from the statements and method of
proof of Theorem~\ref{lem:1.7}. For this, we define
\begin{eqnarray*}
&&\mathcal{M}:=\biggl\{G\dvtx (0,T]\times\mathbb{R}^d\times
\mathbb{R}^d\rightarrow \mathbb{R}; \int_0^t\,ds
\int dz\bigl|G_s(z,y)\beta_{t-s}(x,z)\bigr|<\infty,\\
&&\hspace*{238pt} \forall y\in
\mathbb{R}^d,\forall t\in[0,T]\biggr\}.
\end{eqnarray*}
Furthermore, for $G\in\mathcal{M}$ and $g\in L^\infty([0,T]\times\mathbb
{R}^d)$, we define $Gg(t,x)=\int dy\, g(y)G_t(x,y)$.

\begin{corollary}
\label{cor:38}
Assume that Hypotheses \ref{(h1)} and \ref{(h3)} hold true.
Then $S_a$ is the fundamental solution to the equation $%
f=g+U_{a}f$. That is, for any $g\in L^\infty([0,T]\times\mathbb{R}^d)$,
the solution can be written as
\[
f(t,x)=g(t,x)+\int_{0}^{t}\,ds\int
dy\,g(s,y)S_{a}\bigl((t,x),(s,y)\bigr).
\]
Furthermore, consider the equation
$f=Gg+U_{a}f$ where $g\in L^\infty([0,T]\times\mathbb{R}^d)$
and $G\in\mathcal{M}$.
Then the unique solution can be written as
$
f(t,x)=Gg(t,x)+\int dy\,g(y)\bar S_{a}((t,x),(0,y))
$ where $\bar S_a$ is given by the following uniform absolutely
convergent infinite sum:
%
\begin{eqnarray}
\label{eq:81} \bar S_a\bigl((t,x),(0,y)\bigr)&=&\int
_{0}^{t}\,ds\int dz\, G_s(z,y)S_{a}
\bigl((t,x),(s,z)\bigr)
\nonumber
\\[-8pt]
\\[-8pt]
\nonumber
&=&\sum_{n=1}^\infty
\int_{0}^{t}\,ds\int dz\, G_s(z,y)A_n
\bigl((t,x),(s,z)\bigr).
\end{eqnarray}
\end{corollary}

$\bar S_a$ is usually called the fundamental solution of the equation
$f=Gg+U_{a}f$. We will now start discussing regularity properties of
the solution. We have to replace hypothesis \eqref{eq:h2} by a
slightly stronger hypothesis which will lead to uniform integrability.

\begin{guess}
\label{(h2a')}
%
Given some functions ${\gamma}\dvtx \mathbb{R}_{+}\times\mathbb{R}^d\times
\mathbb{R}^d\rightarrow\mathbb{R}_{+}$ and $G\in\mathcal{M}$. Assume that
there exists $r >0$, $\zeta> 1$ and a
function $\xi\dvtx\mathbb{R}_{+}\rightarrow\mathbb{R}_{+}$ such that the
following holds:

\textup{(i)} For every $z_{0},z_{n}\in\mathbb{R}^d$ and $R>0$ there exists a
constant $C_{R}\equiv C_R(z_0,\break z_n)>0$ such that for every $n\ge2$,
$\delta_{i}>0,i=0,\ldots,n-1$ and $(y_0,y_n)\in B(z_0,r)\times B(z_n,r)$
we have
%
\begin{eqnarray}
\label{H2} &&\int dy_{1}\cdots\int dy_{n-1}1_{\{\sum_{i=1}^{n-1}\llvert
y_{i}\rrvert \leq R\}}
\prod_{i=0}^{n-2}{\gamma}
_{\delta
_{i}}(y_{i},y_{i+1})^{\zeta}\bigl|G_{\delta_{n-1}}(y_{n-1},y_{n})\bigr|^\zeta
\nonumber
\\[-8pt]
\\[-8pt]
\nonumber
&&\qquad\leq C_{R}^{n}\xi \Biggl(\sum
_{i=0}^{n-1}\delta_{i}\Biggr).
\end{eqnarray}

\textup{(ii)} For every $z_{0},z_{n}\in\mathbb{R}^d$ there exists a constant
$C\equiv C(z_0,z_n)>0$ such that for every $n\in\mathbb{N}$, $\delta
_{i}>0,i=0,\ldots,n-1$, $(y_0,y_n)\in B(z_0,r)\times B(z_n,r)$ and
$\varepsilon>0$, there exists $R_\varepsilon>0$ with
%
\begin{eqnarray}
\label{h2a} &&\int dy_{1}\cdots\int dy_{n-1}1_{\{\sum_{i=1}^{n-1}\llvert
y_{i}\rrvert
>R_\varepsilon\}}
\prod_{i=0}^{n-2}{\gamma}
_{\delta
_{i}}(y_{i},y_{i+1})G_{\delta_{n-1}}(y_{n-1},y_{n})
\nonumber
\\[-8pt]
\\[-8pt]
\nonumber
&&\qquad\leq C^n\varepsilon \xi \Biggl(\sum_{i=0}^{n-1}
\delta_{i}\Biggr).
\end{eqnarray}
%
%
\end{guess}

%
%
The reason for both conditions should be clear. The first one, gives a
uniform integrability condition on compact sets. The second condition
states that the measure of the complement of the compact set
$\{\sum_{i=1}^{n-1}\llvert
y_{i}\rrvert \leq R\}$ is sufficiently small.

\begin{lemma}
\label{lemma:1.10}
Assume Hypotheses \ref{(h1)} and \ref{(h3)}.
Suppose that Hypothesis \ref{(h2a')} holds for some $\zeta\in(1,\rho
^{-1})$ and $\gamma$ given in Hypothesis \ref{(h1)}. Furthermore,
assume that $(t,x,y)\rightarrow(G _{t}(x,y), a_t(x,y))$ is continuous
in $(0,T]\times\mathbb{R}^d\times\mathbb{R}^d$. Then
$(t,x,y)\rightarrow\bar S_{a }((t,x),(0,y))$ is continuous.
\end{lemma}

\begin{pf} First recall that \eqref{eq:81} is a uniform absolutely
convergent sum, therefore, it is enough to prove the joint continuity
of each term in the sum. Each term is divided in two integrals on
disjoint sets. The first, on a compact set, is uniformly integrable because
\begin{eqnarray*}
&&\sup_{\llvert  y_{0}\rrvert \llvert  y_{n}\rrvert \leq
K}\int_{s}^{s_{0}}
\,ds_{1} \int dy_{1}\int_{s_1}^{s_{0}}
\,ds_{2}%
 \int dy_{1}\cdots\int_{s_{n-2}}^{s_{0}}
\,ds_{n-1} \int dy_{n-1}%
1_{\{\sum_{i=1}^{n-1}\llvert
y_{i}\rrvert \leq R\}}\\
&&\qquad{}\times \prod
_{i=0}^{n-2}\bigl\llvert A
\bigl((s_{i},y_{i}),(s_{i+1},y_{i+1})
\bigr)\bigr\rrvert ^{\zeta
}\bigl|G_{s_n}(y_{n-1},y_{n})\bigr|^\zeta<
\infty
\end{eqnarray*}
so that we have uniform integrability for the integrand.
Then one may interchange the limit $\lim_{(s_{0},y_{0},y_{n+1})\rightarrow(s_{0}^{\prime},y_{0}^{\prime
},y_{n+1}^{\prime})}$ from outside to inside the integral for fixed
$n$ and $R$.

The argument now finishes fixing $\varepsilon>0$ and, therefore, there
exists $R_\varepsilon$ such that \eqref{h2a} is satisfied. Therefore,
\begin{eqnarray*}
 &&\lim_{\varepsilon\downarrow0}\sup_{(y_{0},y_{n})\in B(z_0,r)\times B(z_{n},r)} \int
_{s}^{s_{0}} \,ds_{1} \int dy_{1}
\int_{s_1}^{s_{0}} \,ds_{2}%
 \int
dy_{1}\cdots\\
&&\qquad{}\times \int_{s_{n-2}}^{s_{0}}
\,ds_{n-1} \int dy_{n-1}%
1_{\{\sum_{i=1}^{n-1}\llvert
y_{i}\rrvert \geq R_\varepsilon\}}\\
&&\qquad{}\times\prod
_{i=0}^{n-2}\bigl\llvert A
\bigl((s_{i},y_{i}),(s_{i+1},y_{i+1})
\bigr) \bigr\rrvert
\bigl|G_{s_n}(y_{n-1},y_{n})\bigr|\le C\lim
_{\varepsilon\downarrow0}\varepsilon\xi(s_0)=0.
\end{eqnarray*}
This gives the continuity of the partial sums and then
of the series itself due to the uniform convergence in \eqref{eq:81}.
\end{pf}

%
%
We discuss now the differentiability properties.

\begin{theorem}
\label{th:1.13}
Assume Hypotheses \ref{(h1)} and \ref{(h3)} and
suppose that\break $\bar G_t(x, y):=\nabla_{y} G_{t}(x,y)$ exists for all
$(t,x,y)\in(0,T]\times\mathbb{R}^d\times
\mathbb{R}^d$. Furthermore, assume that
Hypothesis \ref{(h2a')} is satisfied with $({\gamma},\bar G)$
replacing $(\gamma,G)$ then the application
$y\rightarrow\bar{S}_a((t,x),(0,y))$
is differentiable for $t>0$ and $y\in\mathbb{R}^d$ and the sum below
converges absolutely and uniformly for
$(t,x,y)\in[0,T]\times\mathbb{R}^d\times\mathbb{R}^d$,
\[
\nabla_y \bar{S}_a\bigl((t,x),(0,y)\bigr)=\sum
_{n=1}^\infty \int_{0}^{t}\,ds
\int dz\, \nabla_y G_s(z,y)A_n
\bigl((t,x),(s,z)\bigr).
\]
%
\end{theorem}

The proof is done in a similar way as the proof of Lemma~\ref
{lemma:1.10} using the definition of derivative.

\section{Abstract framework for semigroup expansions}
\label{sec:3}
In this section, we introduce a general framework which will be used in
order to obtain a Taylor-like expansion method for Markovian semigroups.

\begin{guess}
\label{(H1)}
%
$(P_{t})_{t\geq0}$
is a semigroup of linear operators defined on
a space containing $C_{c}^{\infty}(\mathbb{R}^d)$ with infinitesimal generator
$L$ such that $C_{c}^{\infty}(\mathbb{R}^d)\subseteq \operatorname{Dom}(L)$.
$P_{t}f(x)$ is jointly measurable and bounded in the sense that $\llVert
P_{t}f\rrVert  _{\infty}\leq\llVert  f\rrVert  _{\infty}$\mbox{ for
all }%
$f\in C_{c}^{\infty}(\mathbb{R}^d)$\mbox{ and }$t\in{}[0,T]$.
%
\end{guess}

The first goal of this article is to give an expansion for $P_{T}f(x)$
for fixed $%
T>0$ and $f\in C_{c}^{\infty}(\mathbb{R}^d)$ based on a parametrized
semigroup of
linear operators $(P_{t}^{z})_{t\geq0}$, $z\in\mathbb{R}^d$.

In the case of continuous diffusions to be discussed in Section~\ref%
{sec:contdiff}, $P^{z}$ stands for the semigroup of a diffusion process with
coefficients ``frozen'' at $z$. We consider an explicit approximating
class in
the diffusion case in Section~\ref{sec:contdiff} given by the
Euler--Maruyama scheme.

Our hypothesis on $%
(P_{t}^{z})_{t\geq0}$ are:

\begin{guess}
\label{(H2)}
%
For each $z \in\mathbb{R}^d$, $(P_{t}^{z})_{t\geq0}$ is a
semigroup of linear operators defined on a space containing
$C_{c}^{\infty
}(\mathbb{R}^d)$ with infinitesimal generator $L^{z}$ such that
$C_{c}^{\infty}(\mathbb{R}^d)\subseteq \operatorname{Dom}(L^{z})$. We also assume
that $P_{t}^{z}f(x) =\int f(y)p_{t}^{z}(x,y)\,dy$ for any
$f\in C_{c}^{\infty}(\mathbb{R}^d)$, $(x,z)\in\mathbb{R}^d\times\mathbb
{R}^d$ and a jointly measurable probability
kernel $p^{z}\in C((0,T]\times\mathbb{R}^d\times\mathbb{R}^d)$.
%
\end{guess}

The link between $L$ and $L^{z}$ is given by the following hypothesis.

\begin{guess}
\label{(H3)}
$Lf(z)=L^{z}f(z)$ for every $f\in C_c^\infty(\mathbb{R}^d)$ and $z\in
\mathbb{R}^d$.
\end{guess}

To simplify notation, we introduce $Q_{t}f(x):=P_{t}^{x}f(x)$, noticing
that $%
(Q_{t})_{t\geq0}$ is no longer a semigroup but it still satisfies that
$\llVert  Q_{t}f\rrVert  _{\infty}\leq\llVert  f\rrVert  _{\infty}$
for all
$t\in{}[0,T]$. We will use the following notation in the forward
and backward method, respectively
\begin{eqnarray*}
\psi_{t}^{x}(y)&:=&p_{t}^{x}(x,y),\label{eq:psi}
\\
\phi_{t}^{z}(x)&:=&p_{t}^{z}(x,z).\label{eq:phi}
\end{eqnarray*}
The reason for using the above notation is to clarify to which
variables of $p_t^z(x,y)$ an operator applies to. This is the case of,
for example,
$L^z\phi^z_t(x)\equiv(L^z\phi^z_t)(x)$.

The expansion we want to obtain can be achieved in two different ways.
One will be called
the forward method and the other called the backward method. In any of these
methods, the expansion is done based on the semigroup
$(P_{t}^{z})_{t\geq0}$%
, $z\in\mathbb{R}^d$. In the classical Taylor-like expansion one needs
to use polynomials as basis functions. In the forward method, these
polynomials will be replaced by products (or compositions) of the following
basic operator $S$,
\[
S_{t}f(x):=\int\bigl(L^{y}-L^{x}\bigr)f(y)
\psi_t^x(y)\,dy,\qquad f\in\bigcap_{x\in
\mathbb{R}^d}\operatorname{Dom}
\bigl(L^{x}\bigr).
\]
In the backward method, a similar role is played by the operator
%
\begin{equation}
\label{eq:Shat} \hat{S}_{t}f(y):=\int f(x) \bigl(L^x-L^y
\bigr)\phi_t^y(x)\,dx.
\end{equation}

The above
Hypotheses \ref{(H1)}, \ref{(H2)} and \ref{(H3)} will be assumed
throughout
the theoretical part of the article. They will be easily verified in
the examples.

\section{Forward method}
\label{sec:forward}

We first state the assumptions needed in order to implement the forward method.

\begin{guess}
\label{(A0)}
%
$P_{t}^{z}g$, $P_{t}g \in\bigcap_{x\in
\mathbb{R}^d}\operatorname{Dom}(L^{x})$, $\forall g\in C_{c}^{\infty}(\mathbb{R}^d),
z\in\mathbb{R}^d, t\in[0,T]$.
\end{guess}


We assume the following two regularity
properties for the difference operator~$S$.

\begin{guess}
\label{(A3)}
%
There exists a jointly measurable real valued
function $\theta\dvtx (0,T]\times\mathbb{R}^d\times\mathbb{R}^d\rightarrow
\mathbb{R}$, such that for all $f\in C_c^\infty(\mathbb{R}^d)$ we have that
\begin{eqnarray}
S_{t}f(x)=\int f(y)\theta_{t}(x,y)P_{t}^{x}(x,dy)=
\int f(y)\theta_{t}(x,y)p_{t}^{x}(x,y)\,dy,\nonumber\\
 \eqntext{(t,x)
\in(0,T]\times\mathbb{R}^d.}
\end{eqnarray}
We assume that the function $a_t(x,y)=\theta_{t}(x,y)p_{t}^{x}(x,y) $
verifies the Hypotheses~\ref{(h1)} and \ref{(h3)} and that $G_t(x,y)= p^{x}
_{t}(x,y)\in\mathcal{M}$.
%
\end{guess}

Note that the above hypothesis implies that the operator $S$ can be
extended to the space of bounded functions.

\begin{guess}
\label{A3p}
%
For the functions $(a,\gamma)$ and the constant $\rho\in[0,1)$
satisfying Hypothesis \ref{(h1)} and $G_t(x,y)= p^{x}
_{t}(x,y)\in\mathcal{M}$ we assume that the Hypothesis~\ref{(h2a')} is satisfied for some $\zeta\in(1,\rho^{-1})$.
\end{guess}

\begin{remark}
\label{rem:H2}
We remark here that Hypothesis \ref{(A3)} entails some
integration by parts property which will be made clear when dealing with
examples in Section~\ref{sec:contdiff} [see \eqref{def:Sdiff}].
\end{remark}

Define for $(s_{0},x)\in(0,T]\times\mathbb{R}^d$ and $f\in
C_{c}^{\infty}(\mathbb{R}^d)$ the following integral operator:
%
\begin{eqnarray}\qquad
\label{eq:Int}&& I_{s_{0}}^{n}(f) (x)
\nonumber
\\[-8pt]
\\[-8pt]
\nonumber
&&\qquad:=\cases{ %
\displaystyle \int_{0}^{s_{0}}\,ds_{1}
\cdots\int_{0}^{s_{n-1}}\,ds_{n}%
 \Biggl( \prod_{i=0}^{n-1}S_{s_{i}-s_{i+1}}
\Biggr) Q_{s_{n}}f(x), &\quad $\mbox {if } n\geq1,$
\vspace*{2pt}\cr
Q_{s_{0}}f(x), &\quad $\mbox{if } n=0.$}
\end{eqnarray}
We denote by $A_n$ the kernels associated to $a_t(x,y)$ defined in
\eqref{J1}. Then using the change of variables $t_i=s_0-s_i$ we obtain
the following representation
$I^n_{s_0}(f)(x)=\int f(y)I^n_{s_0}(x,y)\,dy$ with
\begin{eqnarray}
\label{eq:11.1} &&{I}_{t_{n+1}}^{n}(x,y)
\nonumber
\\[-8pt]
\\[-8pt]
\nonumber
&&\qquad=\int_{0}^{t_{n+1}}\,dt_{n}
\int p^{y_n} _{t_{n+1}-t_{n}}({y_{n}},y)A_n
\bigl((t_{n+1},x),(t_{n+1}-t_n,y_n)
\bigr)\,dy_{n}.
\end{eqnarray}
%

The following is the main result of this section, which is a
Taylor-like expansion of $P$ based on $Q$.

\begin{theorem}
\label{f1}Suppose that Hypotheses \ref{(A0)} and \ref{(A3)} hold. Then
for every $%
f\in C_{c}^{\infty}(\mathbb{R}^d)$ and $t\in(0,T]$, $I_{t}^{n}(f)$ is
well defined and the sum $%
\sum_{n=1}^{\infty}I_{t}^{n}\times\break (f)(x)$ converges absolutely and uniformly
for $%
(t,x)\in[0,T]\times \mathbb{R}^d$. Moreover,
%
\begin{equation}
P_{t}f(x)=\sum_{n=0}^{\infty}I_{t}^{n}(f)
(x). \label{A7}
\end{equation}
Then for fixed $t\in(0,T]$, $\sum_{n=1}^{\infty}I_{t}^{n}(x,y)$ also
converges absolutely and uniformly for $%
(x,y)\in \mathbb{R}^d\times\mathbb{R}^d$ and we have that
$P_tf(x)=\int f(y)p_t(x,y)\,dy$ where
%
\begin{equation}
\label{A7A}p_{t}(x,y)=p_{t}^{x}(x,y)+\sum
_{n=1}^{\infty}I_{t}^{n}(x,y).
\end{equation}
Furthermore, suppose that $P_tf(x)\ge0$ for $f\ge0$ and $P_t1=1$ for
all $t\ge0$. Then $p_t(x,y)$ is a density function.
\end{theorem}

\begin{pf}
The linear equation on $P_tf$ is obtained, using Hypotheses \ref{(H1)},
\ref{(H2)}, \ref{(A0)} and \ref{(A3)} as follows:
\[
P_{t}f(x)-P_{t}^{x}f(x)=\int
_{0}^{t}\partial_{s_1} \bigl(
P_{t-s_1}^{x}P_{s_1}f \bigr) (x)\,ds_1=
\int_{0}^{t_0}P_{t_0-s_1}^{x}
\bigl(L-L^{x}\bigr)P_{s_1}f(x)\,ds_1.
\]
Note that Hypothesis \ref{(A3)} ensures the finiteness of the above integral.

Using the identity in Hypothesis \ref{(H3)}, $Lg(x)=L^{x}g(x)$ with
$g(y)=P_{s_1}f(y)$, we obtain
\begin{eqnarray*}
P_{t-s_1}^{x}\bigl(L-L^{x}\bigr)P_{s_1}f(x)
&=&\int \bigl(L-L^{x}\bigr)P_{s_1}f(y)P_{t-s_1}^{x}(x,dy)\\
&=&
\int \bigl(L^{y}-L^{x}\bigr)P_{s_1}f(y)P_{t-s_1}^{x}(x,dy)
\\
&=&S_{t-s_1}P_{s_1}f(x).
\end{eqnarray*}
Therefore, we have the following equation:
%
\begin{equation}
\label{eq:11.1a} P_{t}f(x)=P_{t}^{x}f(x)+\int
_0^{t}\,ds_1 \int dy\,
P_{s_1}f(y)\theta_{t-s_1}(x,y)p_{t-s_1}^{x}(x,y).
\end{equation}
This is equation \eqref{eq:app1} with $a_t(x,y)=\theta
_{t}(x,y)p_{t}^{x}(x,y)$. Therefore, due to Hypotheses \ref{(A0)} and
\ref{(A3)} we obtain that the hypotheses needed for the application of
Lemma~\ref{lemma:1.5} and Theorem~\ref{lem:1.7} are satisfied.
Therefore, we obtain that
$\sum_{n=0}^{\infty}I_{t}^{n}(f)(x)$ converges absolutely and
uniformly and is the unique solution
of~\eqref{eq:11.1a}.

Corollary~\ref{cor:38} gives (\ref{A7A}).
Finally, one proves that as
the semigroup $P$ is positive then $p_{{t}}(x,y)$ has to be positive
locally in $y$ for fixed $({t},x)$ then as $P_{{t}}1=1$ one obtains that
$\int p_{{t}}(x,y)\,dy=1$.
\end{pf}

\subsection{Probabilistic representation using the forward
method}
\label{sec:probaforw}

Our aim now is to give a probabilistic representation for the formula
(\ref%
{A7}) that may be useful for simulation.


\begin{guess}
\label{(HMC)}
There exists a continuous Markov process
$X^\pi=\break \{X^\pi_t;t\in[0,T]\}$ such that $X_{0}^{\pi}=x$ and for any
$t>s$ %
\begin{equation}
P\bigl(X_{t}^{\pi}\in dy'\mid
X_{s}^{\pi
}=y\bigr)=P_{t-s}^{y}
\bigl(y,dy'\bigr)=p_{t-s}^y
\bigl(y,y'\bigr)\,dy'. \label{A12}
\end{equation}
\end{guess}

With this assumption, we have that $S_{1}f(x)=\mathbb{E} [
f(X_{t_{1}}^{\pi
})\theta_{t_{1}}(x,X_{t_{1}}^{\pi}) ] $ and $%
Q_{t_{1}}f(x)=\mathbb{E}[f(X_{t_{1}}^{\pi})]$. Therefore, using these
representations, we obtain the probabilistic representation of the integrand
in (\ref{eq:Int}):
\[
\Biggl( \prod_{j=0}^{n-1}S_{t_{j+1}-t_{j}}
\Biggr) Q_{T-t_{n}}f(x)=\mathbb {E} \bigl[ f\bigl(X_{T}^{\pi}
\bigr)\theta_{t_{n}-t_{n-1}}\bigl(X_{t_{n-1}}^{\pi},X_{t_{n}}^{\pi
}
\bigr)\cdots\theta_{t_{1}-t_{0}}\bigl(X_{t_{0}}^{\pi},X_{t_{1}}^{\pi}
\bigr) \bigr].
\]

Finally, to obtain the probabilistic interpretation for the representation
formula~(\ref{A7}), we need to find the probabilistic representation of the
multiple integrals in~(\ref{eq:Int}).

For this, we consider a Poisson
process $(J_{t})_{t\ge0}$ of parameter $\lambda=1$ and we denote by
$\tau_{j}$,
$j\in\mathbb{N}$, its jump times (with the convention that $\tau_{0}=0)$.
%
Conditionally to $J_{T}=n$, the jump times are distributed as the order
statistics of a sequence $n$ independent uniformly distributed random
variables on $[0,T]$. Therefore, the multiple integrals in (\ref{eq:Int})
can be interpreted as the expectation taken with respect to these jump times
given that $J_{T}=n$. Therefore, for $n\geq1$ we have
\[
I_{T}^{n}(f) (x)=e^{T}\mathbb{E} \Biggl[
1_{\{J_{T}=n\}}f\bigl(X_{T}^{\pi
}\bigr)\prod
_{j=0}^{n-1}\theta_{\tau_{j+1}-\tau_{j}}
\bigl(X_{\tau_{j}}^{\pi
},X_{\tau_{j+1}}^{\pi}\bigr)
\Biggr],
\]
where (with a slight abuse of notation), $\pi$ denotes the random time
partition of $[0,T]$, $\pi\equiv\pi(\omega)= \{ \tau_{i}(\omega
)\wedge T;i=0,\ldots,J_{T}(\omega)+1 \} $.

From now on, in order to simplify the notation, we denote $\tau
_{T}\equiv
\tau_{J_{T}}$. Given the above
discussion, we have the main result for this section.

\begin{theorem}
\label{th:forwardrep}Suppose that Hypotheses \ref{(A0)}, \ref{(A3)} and
\ref{(HMC)} hold. Recall that $\psi_t^x(y)=p_t^x(x,y)$ and define
$\Gamma_{T}(x)\equiv\Gamma_{T}(x)(\omega)$ as
\[
\Gamma_{T}(x)=\cases{ %
\displaystyle\prod
_{j=0}^{J_{T}-1}\theta_{\tau_{j+1}-\tau_{j}}
\bigl(X_{\tau_{j}}^{\pi
},X_{\tau_{j+1}}^{\pi}\bigr), &\quad
$\mbox{if } J_{T}\geq1,$
\vspace*{2pt}\cr
1, &\quad $\mbox{if } J_{T}=0.$}
\]
Then the
following probabilistic representations are satisfied for $f\in
C_{c}^{\infty
}(\mathbb{R}^d)$:
%
\begin{eqnarray}
P_{T}f(x)&=&e^{T}\mathbb{E} \bigl[ f\bigl(X_{T}^{\pi}
\bigr)\Gamma_{T}(x) \bigr], \label{A13}
\\
p_{T}(x,y)&=&e^{T}\mathbb{E} \bigl[ \psi_{T-\tau_{T}}^{X_{\tau
_{T}}^{\pi}}(y)
\Gamma _{T}(x) \bigr]. \label{A10}
\end{eqnarray}
\end{theorem}

\begin{remark}
%
%
1. Extensions for bounded measurable functions $f$ can be obtained if
limits are taken in \eqref{A13}.

2. The above representations (\ref{A10}) and (\ref{A13}) may be
obtained using
a Poisson process of arbitrary parameter $\lambda>0$ instead of
$\lambda
=1$. In fact, if we denote $\{J_{t}^{\lambda}, t\ge0\}$ a Poisson
process, by $\tau
_{i}^{\lambda}$ the jump times and by $\pi_{\lambda}$ the corresponding
random time grid. Then the formula (\ref{A10}) becomes
\[
p_{T}(x,y)=e^{\lambda T}\mathbb{E} \bigl[ \lambda^{-J_{T}^{\lambda
}}
\psi_{T-\tau
_{J_{T}^{\lambda}}}^{X_{\tau_{T}}^{\pi_{\lambda}}}(y)\Gamma_{T}(x)%
 \bigr].
\label{A10'}
\]
\end{remark}

\subsection{Regularity of the density using the forward method}
\label{sec:rgul}
Now that we have obtained the stochastic representation, we will
discuss the differentiability of
$p_T(x,y)$ with respect to $y$. This type of property is also proved
when the analytical version of the parametrix method is discussed in
the particular case of fundamental solutions of parabolic PDEs (see,
e.g., Chapter~1 in \cite{F2}).

\begin{theorem}
\label{f3}Suppose that the Hypotheses \ref{(A0)}, \ref{(A3)}, \ref{A3p}
are satisfied. Furthermore assume that $(t,x,y)\rightarrow
(p^x_t(x,y),a_t(x,y))$ is continuous in\break $(0,T]\times\mathbb{R}^d\times
\mathbb{R}^d$.
Then $(t,x,y)\rightarrow p_{t}(x,y)$ is continuous on $(0,T]\times
\mathbb{R}^d\times\mathbb{R}^d$.
\end{theorem}

\begin{theorem}
\label{f4}
Suppose that the Hypotheses \ref{(A0)}, \ref{(A3)} and \ref{A3p} are
satisfied. Furthermore, we assume that $y\rightarrow p_t^x(x,y)$ is
differentiable for all $(t,x)\in\mathbb{R}_+\times\mathbb{R}^d$ and
that Hypothesis \ref{(h2a')} is satisfied for $\nabla_yp_t^x(x,y)$
instead of $G$.
Then for every $(t,x)\in(0,T]\times\mathbb{R}^d$, the function
$y\rightarrow p_{t}(x,y)$ is differentiable. Moreover,
\[
\nabla_{y}p_{T}(x,y)=e^{T}\mathbb{E} \bigl[
\nabla_{y}p _{T-\tau
_{T}}^{X_{\tau
_{T}}^{\pi}}\bigl(X^\pi_{\tau
_{T}},y
\bigr)\Gamma_{T}(x) \bigr].
\]
\end{theorem}
%
\section{The backward method: Probabilistic representation using the
adjoint semigroup}
\label{sec:back}

We will now solve equation \eqref{eq:3.1} in dual form. We start with a
remark in order to
beware the reader about the nonapplicability of the forward method
directly to the dual problem.

Usually, in semigroup theory one assumes that for each $t>0$, $P_{t}$
maps continuously $%
L^{2}(\mathbb{R}^d)$ into itself.
Then $P_{t}^{\ast}$ can be defined and it is still a semigroup which
has as infinitesimal
operator $L^{\ast}$ defined by $ \langle L^{\ast}g,f \rangle
= \langle g,Lf \rangle$ for $f$, $g\in C_{c}^{\infty}(\mathbb
{R}^d)$. Assume, for the sake of the present discussion, that for every
$x\in\mathbb{R}^d$, $P_{t}^{x}$ maps continuously $%
L^{2}(\mathbb{R}^d)$ into itself and we define $P_{t}^{x,\ast}\equiv
 (
P_{t}^{x} ) ^{\ast}$ and $L^{x,\ast}$ by $ \langle
P_{t}^{x,\ast
}g,f \rangle= \langle g,P_{t}^{x}f \rangle$ and $%
 \langle L^{x,\ast}g,f \rangle= \langle g,L^{x}f
\rangle
$ for $f$, $g\in C_{c}^{\infty}(\mathbb{R}^d)$.

Our aim is to obtain for $P^{\ast}$ a representation which is similar to
the one obtained for $P$ in Theorem~\ref{th:forwardrep}$.$ Unfortunately,
the adjoint version of the arguments given in Section~\ref{sec:forward} do
not work directly. In fact, if $P_{t}^{x,\ast}$ denotes the adjoint
operator of $%
P_{t}^{x}$ then the relation $Lg(x)=L^{x}g(x)$ does not imply $L^{\ast
}g(x)= ( L^{x} ) ^{\ast}g(x)$. To make this point clearer, take,
for example, the case of a one-dimensional diffusion process with
infinitesimal operator $Lg(y)=a(y)\triangle g(y)$ then $L^{x}g(y)=a(x)%
\triangle g(y)$ (for more details, see Section~\ref{sec:contdiff}).
Then $%
L^{\ast}g(y)=\triangle(ag)(y)$ and $ ( L^{x} ) ^{\ast
}g(y)=a(x)\triangle g(y)=L^{x}g(y)$. So, letting the coefficients of $%
L^{\ast}$ be frozen at $y=x$ does not coincide with $ ( L^{x} )
^{\ast}$ and, therefore, the previous argument will fail.\footnote{%
Note that if we wanted to freeze coefficients as in the forward method one
may be lead to the study of the operator $L^{\ast,z}g(y)=a(z)\Delta
g(y)+2 \langle\nabla a(z),\nabla g(y) \rangle+g(y)\Delta a(z)$.
Although this may have an interest in itself, we do not pursue this
discussion here as this
will again involve derivatives of the coefficients while in this
section we
are pursuing a method which may be applied when the coefficients are H\"
older continuous.} In order not to
confuse the reader, we will keep using the superscript $\ast$ to
denote adjoint operators
while other approximating operators will be denoted by the superscript
$\hat{}$ (hat).


Note that in the diffusion case, proving that for each $t>0$, $P_{t}$
maps continuously $%
L^{2}(\mathbb{R}^d)$ into itself is not easy in general. Therefore,
instead of adding this as a hypothesis, we will make additional
hypotheses on the approximation process. This issue will demand us to
introduce hypotheses that did not have any counterpart in the forward
method. Still, once the linear Volterra equation is obtained, the
arguments are parallel and we will again use the results in Section~\ref{sec:FE}. Let us introduce some notation and the main hypotheses. We
define the linear operator
\[
\hat{Q}_{t}f(y):= \bigl( P_{t}^{y} \bigr)
^{\ast}f(y)=\int f(x)p_{t}^{y}(x,y)\,dx=\int f(x)
\phi_{t}^{y}(x)\,dx.
\]
We assume the following.

\begin{guess}
\label{A0ast}
%
\textup{(0)} $P_tf(x)=\int f(y)P_t(x,dy)$ for all $f\in C_b(\mathbb{R}^d)$. In
particular, $P_t$ is an integral operator.\vspace*{-6pt}
\begin{longlist}[(iii)]
\item[(i)]$\int p_{t}^{y}(x,y)\,dy<\infty$, for
all $x\in\mathbb{R}^d$ and $\int p_{t}^{y}(x,y)\,dx<\infty$ for all
$y\in
\mathbb{R}^d$.

\item[(ii)] $\lim_{\varepsilon\rightarrow0}P_{T+\varepsilon}^{z,\ast}g(w)
=P_{T}^{z,\ast}g(w)$ and
$\lim_{\varepsilon\rightarrow0}\int h(z)\phi^z_\varepsilon
(w)\,dz=h(w)$ for all $(z,w)\in\mathbb{R}^d\times\mathbb{R}^d$ and for
$g, h\in C_{c}^\infty(\mathbb{R}^d)$.

\item[(iii)] $\phi_{t}^{z} \in \operatorname{Dom}(L)\cap(\bigcap_{y\in\mathbb
{R}^d}\operatorname{Dom}(L^{y}))$, for all $(t,z)\in(0,T]\times
\mathbb{R}^d$.
%
\end{longlist}
\end{guess}

With these definitions and hypotheses, we have by the semigroup
property of $%
P^{z}$, that
%
\begin{equation}
P_{t}^{z}\phi_{\varepsilon}^{z}(x)=p_{t+\varepsilon}^{z}(x,z)=
\phi _{t+\varepsilon}^{z}(x). \label{Ad1}
\end{equation}
As stated before, we remark that $\hat{Q}\neq Q^{\ast}$. In fact, $\hat{Q}$
is defined through a density whose coefficients are ``frozen'' at the
arrival point of the underlying process. Also note that due to
Hypothesis \ref{A0ast} then $\llVert  \hat{Q}_{t}f\rrVert  _{\infty
}\leq C_{T}\llVert  f\rrVert
_{\infty}$ for all $t\in{}[0,T]$.

Before introducing the next two hypotheses, we explain the reasoning
behind the notation to follow. In the forward method, it was clear that
the dynamical system expressed through the transition densities went
from a departure point $x$ to an arrival point $y$ with transition
points $y_i$, $i=0,\ldots,n+1$, $y_0=x$ and $y_{n+1}=y$. In the backward
method, the situation is reversed. The initial point for the method is~$y$,
the arrival point is $x$ and $y_0=y$ and $y_{n+1}=x$. The notation
to follow tries to give this intuition.

\begin{guess}
\label{A1ast}
%
We suppose that there exists a continuous
function $\hat{\theta} \in C((0,T]\times\mathbb{R}^d\times\mathbb{R}^d)$
such that $(L^{y_1}-L^{y})\phi_{t}^{y}(y_1) =\hat{\theta
}_{t}(y_1,y)\phi
_{t}^{y}(y_1)$. Moreover, we assume that $\hat{\theta}_{t}(y_1,y)\phi
_{t}^{y}(y_1)$ is integrable $P_{s}(x,dy_1)$ for all $(s,x)\in
[0,T]\times\mathbb{R}^d$ and $(t,y)\in(0,T]\times\mathbb{R}^d$.
%
\end{guess}

Define the function ${\hat{a}}_{t}(x,y):=\hat{\theta
}_{t}(y,x)p_{t}^x(y,x)=\hat{\theta}_{t}(y,x)\phi_t^x(y)$.

\begin{guess}
\label{A2ast}
%
Assume that the function ${\hat{a}}$ satisfies the Hypotheses
\ref{(h1)} and~\ref{(h3)} with $G_t(x,y)=p^x(y,x)\in\mathcal{M}$.
Furthermore, we assume that the corresponding function $\gamma$ satisfies
$\sup_{(t,y)\in[0,T]\times\mathbb{R}^d}\int|\gamma_t(x,y)| \,dx<\infty$
and that there exists $\zeta\in(1,\rho^{-1})$ such that for every $R>0$
%
\begin{equation}
\label{h2b} \sup_{(t,y)\in[0,T]\times\mathbb{R}^d}\int1_{\{|x|\le R\}}\bigl|\gamma
_t(x,y)\bigr|^\zeta \,dx<\infty.
\end{equation}
%
%
\end{guess}

\begin{guess}
\label{A2astA}
%
For the function $ {\hat{a}}_t(x,y)$ we assume that Hypothesis
\ref{(h2a')} is satisfied for some $\zeta\in(1,\rho^{-1})$.
%
\end{guess}

We define now [recall \eqref{eq:Shat} and Hypothesis \ref{A1ast}]
\[
\hat{S}_{t}f(y):=\int f(x)\hat{a}_{t}(y,x)\,dx.
\label{Ad4}
\]
%
For $g\in C_{c}^{\infty}(\mathbb{R}^d)$, we define
%
\begin{eqnarray}\label{Aa1}
&&\hat{I}_{s_{0}}^{n}(g) (y)
\nonumber
\\[-8pt]
\\[-8pt]
\nonumber
&&\qquad:=\cases{ %
\displaystyle \int_{0}^{s_{0}}\,ds_{1}\cdots\int
_{0}^{s_{n-1}}\,ds_{n}%
 \Biggl( \prod
_{i=0}^{n-1}\hat{S}_{s_{i}-s_{i+1}} \Biggr)
\hat {Q}_{s_{n}}g(y), &\quad $\mbox{if } n\geq1,$
\vspace*{2pt}\cr
\hat{Q}_{s_{0}}g(y), &\quad $\mbox{if }n=0.$}
\end{eqnarray}
Furthermore, we define the adjoint operators
\begin{eqnarray*}
\hat{Q}_{t}^{\ast}f(x) &:=&\int f(y)p_{t}^{y}(x,y)\,dy,
\\
\hat{S}_{t}^{\ast}f(x) &:=&\int f(y)\hat{a}_{t}(y,x)\,dy.
\end{eqnarray*}
Note that due to the Hypotheses \ref{A0ast}(i) and \ref{A2ast} we have
that for any $f\in L^{\infty}$
%
\begin{eqnarray}
\sup_{t}\bigl\llVert \hat{Q}_{t}^{\ast}f
\bigr\rrVert _{\infty} &\leq&C\llVert f\rrVert _{\infty}
,\label{eq:qhats}
\\
\bigl\llVert \hat{S}_{t}^{\ast}f\bigr\rrVert _{\infty}
&\leq&\frac{C}{t^{\rho
}}%
\llVert f\rrVert _{\infty}.\label{eq:shats}
\end{eqnarray}
As in \eqref{eq:Int}, we define the following auxiliary operators for
\begin{eqnarray*}
&&\hat{I}^{n,*}_{t}(f):=\cases{ %
\displaystyle\int_{0}^{t}\,dt_{n}\cdots\int
_{0}^{t_{2}}\,dt_{1}\hat{Q}_{t_{1}}^{\ast}%
\hat{S}_{t_{2}-t_{1}}^{\ast}\cdots\hat{S}_{t-t_{n}}^{\ast}f,
& \quad $n\geq1,$
\vspace*{2pt}\cr
\hat{Q}_{t}^{\ast}f, &\quad  $n=0,$ }
\\
%
&&\hat{I}_{t}^{n}(y_{0},y_{n+1})\\
&&\qquad:=\int
_{0}^{t}\,dt_{n}\cdots%
\int
_{0}^{t_{2}}\,dt_{1}\int dy_{1}
\cdots\int dy_{n}\,{\hat{A}_n}\bigl((t,y_0),(t-t_n,y_n)
\bigr)p_{t_n}^{y_n}(y_{n+1},y_{n}).
\end{eqnarray*}
Here, ${\hat{A}_n}$ denotes the same function defined in \eqref{J1}
where ${\hat{a}}$ is used instead of $a$.
Note that $\langle\hat{I}^{n,*}_{t_{0}}(f),g\rangle=\langle f,\hat
{I}^n_{t_0}g\rangle$ and $ \hat{I}^{n,*}_t(f)(x)=\int f(y)\hat
{I}^{n,*}_{t}(y,x)\,dy$ for $f, g\in C_c^\infty(\mathbb{R}^d)$. Our main
result in this section is:

\begin{theorem}
\label{b1}Suppose that Hypotheses \ref{A0ast}, \ref{A1ast} and \ref
{A2ast}. Then for every $g\in
C_c^{\infty}(\mathbb{R}^d)$ the sum $\sum_{n=0}^{\infty}\hat{I}_{t}^{n}(g)(y)$
converges absolutely and uniformly for $(t,y)\in(0,T]\times\mathbb
{R}^d$ and
the following representation formula is satisfied:
%
\begin{equation}
P_{t}^{\ast}g(y)=\sum_{n=0}^{\infty}
\hat{I}_{t}^{n}(g) (y),\qquad dy\mbox{-a.s.}, t\in(0,T].
\label{Ad9f}
\end{equation}
The above equality is understood in the following weak sense $
\langle
P_{t}^{\ast}g,h \rangle= \langle g,P_{t}h \rangle= \sum_{n=0}^{\infty} \langle\hat{I}_{t}^{n}(g),h \rangle$
for all $(g,h)\in C_c^\infty(\mathbb{R}^d)\times C_{c}^\infty(\mathbb{R}^d)$.
Furthermore, $\sum_{n=0}^{\infty}\hat{I}^{n,*}_{t}(f)(x)$ converges
absolutely and uniformly for $x\in\mathbb{R}^d$\vspace*{1pt} and fixed $f\in
C_c^\infty(\mathbb{R}^d)$, $t\in(0,T]$ and it satisfies
\[
P_{t}f(x)=\sum_{n=0}^{\infty}
\hat{I}^{n,*}_t(f) (x),  \qquad dx\mbox{-a.s.}
\]
Finally, $\sum_{n=0}^{\infty}\hat{I}^{n}_{t}(y,x)$ converges
absolutely and uniformly for $(x,y)\in\mathbb{R}^d\times\mathbb{R}^d$
for fixed $t>0$ and there exists a jointly measurable function
$p_t(x,y)$ such that we have that for $f\in C_c^\infty(\mathbb{R}^d)$
we have $P_tf(x)=\int f(y)p_t(x,y)\,dy$ and it is given by
\[
p_t(x,y)=p_t^y(x,y)+\sum
_{n=1}^\infty\hat{I}^n_t(y,x).
\]
Furthermore, suppose that $P_tf(x)\ge0$ for $f\ge0$ and $P_t1=1$ for
all $t\ge0$. Then $p_t(x,y)$ is a density function.
\end{theorem}

\begin{pf} Many of the arguments are similar to the proof of
Theorem~\ref{f1}.
In fact, we first establish the Volterra equations satisfied by
$P_t^*$. In
order to do this, we need an approximation argument. We fix
$\varepsilon>0$
and we recall that due to Hypotheses \ref{(H2)} and \ref{A0ast}(iii),
we have
for each $z\in\mathbb{R}^d$ that $P_{T-s}^{z}\phi_{\varepsilon
}^{z}=\phi_{T-s+\varepsilon
}^{z}=p_{s+\varepsilon}^{z}(\cdot,z)\in \operatorname{Dom}(L)$ and $Lf(x)=L^{x}f(x)$. Then
from Hypotheses \ref{A1ast} and \ref{A2ast}, we have that for $0<s<T$,
\begin{eqnarray*}
\partial_{s} \bigl( P_{s}P_{T-s}^{y_0}
\bigr) \phi_{\varepsilon}^{y_0}(x) &=&P_{s}
\bigl(L-L^{y_1}\bigr)P_{T-s}^{y_0}
\phi_{\varepsilon}^{y_0}(x)\\
&=&\int P_{s}(x,dy_1)
\bigl(L-L^{y_0}\bigr)P_{T-s}^{y_0}
\phi_{\varepsilon}^{y_0}(y_1)
\\
&=&\int P_{s}(x,dy_1)\hat{a}_{T-s+\varepsilon}(y_0,y_1).
\end{eqnarray*}
We take $g, h\in C_c^{\infty}(\mathbb{R}^d)$ and we note that due to
Hypothesis \ref{A2ast}
\begin{eqnarray*}\label{eq:adjnormest}
&&\int dx\bigl\llvert g(x)\bigr\rrvert \int P_{s}(x,dy_1)
\int dy_0\bigl\llvert h(y_0)\bigr\rrvert \bigl\llvert
\hat{a}_{T-s+\varepsilon}(y_0,y_1)\bigr\rrvert \\
&&\qquad\leq
C_{3}(T-s+\varepsilon)^{-\rho}\llVert h\rrVert _{\infty}
\int dx\bigl\llvert g(x)\bigr\rrvert \int P_{s}(x,dy_1)
\\
&&\qquad \leq C_{3}(T-s)^{-\rho}\llVert h\rrVert _{\infty}
\llVert g\rrVert _{1}.
\nonumber
\end{eqnarray*}
The above expression is integrable with respect to $1_{(0,T)}(s)\,ds$ for
$%
\rho\in(0,1)$. Therefore this ensures that Fubini--Tonelli's theorem
can be
applied and multiple integrals appearing in any order will be well defined.\vadjust{\goodbreak}

Furthermore, by Hypotheses \ref{A1ast}, \ref{A2ast} [see \eqref{h2b}]
and the fact that $h\in C_{c}^\infty(\mathbb{R}^d)$, we have that for
fixed $s\in{}[0,T)$
we can take limits as $\varepsilon\rightarrow0$ for $\int dy_0\llvert
h(y_0)\rrvert\times\break   \llvert
\hat{a}_{T-s+\varepsilon}(y_0,y_1)\rrvert $, and that the uniform
integrability property is satisfied. Therefore, we finally obtain that the
following limit exists, is finite and the integration order can be exchanged
so that
\begin{eqnarray*}
&&\lim_{\varepsilon\rightarrow0}\int_{0}^{T} \,dt
\int dy_0\,h(y_0)\int dx\,g(x)\int P_{s}(x,dy_1)
\hat{a}_{T-s+\varepsilon}(y_0,y_1)\\
&&\qquad=\int
_{0}^{T} \,dt\int dy_0\,h(y_0)
\int dx\,g(x)\int P_{s}(x,dy_1)\hat{a}%
_{T-s}(y_0,y_1).
\end{eqnarray*}
From the previous argument, the following sequence of equalities are valid
and the limit of the right-hand side below exists:
%
\begin{eqnarray}
\label{eq:appexp1}&& \int dy_0\,h(y_0) \bigl( \bigl\langle
g,P_{T}\phi_{\varepsilon}^{y_0} \bigr\rangle - \bigl\langle
g,P_{T}^{y_0}\phi_{\varepsilon}^{y_0} \bigr\rangle
\bigr)\nonumber\\
 &&\qquad=\int dy_0\,h(y_0)\int dx\,g(x)\int
_{0}^{T}\partial_{t} \bigl(
P_{s}P_{T-s}^{y_0} \bigr) \phi_{\varepsilon}^{y_0}(x)\,dt
\\
&&\qquad=\int_{0}^{T}\,dt\int dy_0\,h(y_0)
\int dx\,g(x)\int P_{s}(x,dy_1)\hat{a}%
_{T-s+\varepsilon}(y_0,y_1).
\nonumber
\end{eqnarray}
In order to obtain the linear Volterra type equation, we need to take
limits in (%
\ref{eq:appexp1}). To deal with the limit of the left-hand side of (\ref
{eq:appexp1}), we note that given the assumptions $g, h\in C_c^\infty
(\mathbb{R}^d)$ and Hypothesis \ref{A0ast}(ii), we have
\begin{eqnarray*}
\lim_{\varepsilon\rightarrow0}\int dy_0\,h(y_0) \bigl
\langle g,P_{T}\phi _{\varepsilon}^{y_0} \bigr\rangle&=&\lim
_{\varepsilon\rightarrow
0}\int g(y_1)\int P_{T}(y_1,dw)
\int dy_0\,h(y_0)\phi_{\varepsilon
}^{y_0}(w)\\
&=&
\langle P_{T}h,g \rangle,
\\
\lim_{\varepsilon\rightarrow0}\int dy_0\,h(y_0) \bigl
\langle P_{T}^{y_0,
\ast
}g,\phi_{\varepsilon}^{y_0} \bigr
\rangle&=&\lim_{\varepsilon
\rightarrow
0}\int dy_0\,h(y_0)P_{T+\varepsilon}^{y_0, \ast}g(y_0)\\
&=&
\int dy_0\,h(y_0)P_{T}^{y_0,\ast}g(y_0).
\end{eqnarray*}
Therefore, taking limits in (\ref{eq:appexp1}), we obtain
\begin{eqnarray*}
\langle P_{T}h,g \rangle&=&\int dy_0\,h(y_0)P_{T}^{\ast
,y_0}g(y_0)
\\
&&{}+\int dy_0\,h(y_0)\int dx\,g(x)\int_{0}^{T}\,ds
\int P_{s}(x,dy_1)\hat{a}%
_{T-s}(y_0,y_1)
\\
&=& \bigl\langle\hat{Q}_{T}^{\ast}h,g \bigr\rangle +\int
_{0}^{T} \bigl\langle P_{s}
\hat{S}_{T-s}^{\ast}h,g \bigr\rangle \,ds.
\end{eqnarray*}
Rewriting this equation with the adjoint of a densely defined operator,
we obtain the Volterra-type equation
\[
P_T^*g(x)=\hat{Q}_Tg(x)+\int_0^T\,ds
\int dyP_s^*g(y){\hat{a}}_{T-s}(x,y).
\]
This equation has a solution due to the results in Section~\ref{sec:FE}
as we have made the necessary hypotheses to apply the results of
Corollary~\ref{cor:38}.
Therefore, it follows that (\ref{Ad9f}) is the unique solution of the
above equation.
The proof of the other statements are done in the same way as in the
proof of Theorem~\ref{f1}.
\end{pf}

\begin{remark}
The previous proof is also valid with weaker conditions on $g$ and $h$.
For example, $g\in L^1(\mathbb{R}^d)\cap L^\infty(\mathbb{R}^d)$ and
$h\in C_b (\mathbb{R}^d)$ will suffice with an appropriate change of hypothesis.
\end{remark}

\subsection{Probabilistic representation and regularity using the
backward method}
\label{sec:5.1}
We deal now with the representation of the density associated with the
semigroup~$P_{T}$. We recall that in the Section~\ref{sec:probaforw}
[see (\ref{A12})] we have performed a similar construction.

\begin{guess}
\label{(HMCB)}
There exists a continuous Markov process
$\{X_{t}^{*, \pi}(y),t\in[0,T]\}$, $y\in\mathbb{R}^d$ such that
$X_{0}^{*, \pi}(y)=y$ and for any $t>s$ we have
\begin{eqnarray*}
P\bigl(X_{t}^{\ast,\pi} (y)\in dy_2\mid
X_{s}^{\ast,\pi
}(y)=y_1\bigr)&=&C_{t-s}^{-1}(y_1)P_{t-s}^{y_1,\ast
}(y_1,dy_2)\\
&=&C_{t-s}^{-1}(y_1)
\phi _{t-s}^{y_1}(y_2)\,dy_2,
\\
C_{t-s}(y_1) &:=&\int\phi_{t-s}^{y_1}(y_2)\,dy_2.
\end{eqnarray*}
\end{guess}
Let $(J_{t})_{t\ge0}$ be a Poisson process of parameter $\lambda=1$
and we denote by $\tau_{j}$,
$j\in\mathbb{N}$, its jump times (with the convention that $\tau_{0}=0)$.
Then the same
arguments as in the previous section give the representation
\begin{eqnarray*}
&&I_{T}^{\ast,n}(g) (y)\\
&&\qquad=e^{T}\mathbb{E} \Biggl[
1_{\{J_{T}=n\}
}g\bigl(X_{T}^{\ast,\pi
}(y)\bigr)C_{T-\tau_{n}}
\bigl(X_{\tau_{n}}^{\ast,\pi}(y)\bigr)\\
&&\hspace*{24pt}\qquad\quad{}\times\prod_{j=0}^{n-1}C_{\tau
_{j+1}-\tau_{j}}
\bigl(X_{\tau_{j}}^{\ast,\pi}(y)\bigr)\hat{\theta}_{\tau
_{j+1}-\tau
_{j}}
\bigl(X_{\tau_{j+1}}^{\ast,\pi}(y),X_{\tau_{j}}^{\ast,\pi
}(y)
\bigr) \Biggr].
\end{eqnarray*}
We define
\begin{eqnarray*}
&&\Gamma_{T}^{\ast}(y)\\
&&\qquad=\cases{ %
\displaystyle
C_{T-\tau_{J_{T}}}\bigl(X_{\tau_{J_{T}}}^{\ast,\pi
}(y)\bigr)\prod
_{j=0}^{J_{T}-1}C_{\tau_{j+1}-\tau_{j}}\bigl(X_{\tau_{j}}^{\ast
,\pi}(y)
\bigr)%
\hat{\theta}_{\tau_{j+1}-\tau_{j}}\bigl(X_{\tau_{j+1}}^{\ast,\pi
}(y),X_{\tau
_{j}}^{\ast,\pi}(y)
\bigr), \vspace*{2pt}\cr
\qquad \mbox{if } J_{T}\geq1,
\vspace*{2pt}\cr
C_{T}(y), \vspace*{2pt}\cr
\qquad \mbox{if } J_{T}=0.}
\end{eqnarray*}
Sometimes we may use the notation $X_{\tau_{j}}^{\ast,\pi}(y)$ to
indicate that $X_{0}^{\ast,\pi}(y)=y$. The main result in this
section is
about representations of the adjoint semigroup $P^{\ast}$ and its densities.

\begin{theorem}
Suppose that Hypotheses \ref{A0ast}, \ref{A1ast} and \ref{(HMCB)} hold
then the following
representation formula is valid for any $g\in C_{c}^{\infty}(\mathbb{R}^d)$:
%
\begin{eqnarray} \label{Ad6}
P_{T}^{\ast}g(y)&=&P_{T}^{\ast,z}g(y)+e^{T}
\mathbb{E} \bigl[ g\bigl(X_{T}^{\ast,\pi
}(y)\bigr)
\Gamma_{T}^{\ast}(y)1_{\{J_{T}\geq1\}} \bigr]
\nonumber
\\[-8pt]
\\[-8pt]
\nonumber
& =&e^{T}
\mathbb {E} \bigl[ g\bigl(X_{T}^{\ast,\pi}(y)\bigr)
\Gamma_{T}^{\ast}(y) \bigr].
\end{eqnarray}
\end{theorem}

\begin{theorem}
\label{th:3.4}Suppose that Hypotheses \ref{A0ast}, \ref{A1ast} and \ref{(HMCB)}
hold then the
following representation formula for the density is valid:
%
\begin{equation}
p_{T}(x,y)=e^{T}\mathbb{E} \bigl[ p _{T-\tau_{T}}^{X_{\tau_{T}}^{\ast
,\pi
}(y)}
\bigl(x,X_{\tau_{T}}^{\ast,\pi
}(y)\bigr)\Gamma_{T}^{\ast}(y)
\bigr]. \label{eq:adjden}
\end{equation}
In particular, let $Z$ be a random variable with density $h\in
L^{1}(\mathbb{R}^d;\mathbb{R}_+)$ then we have
\[
P_{T}h(x)=e^{T}\mathbb{E} \bigl[ p _{T-\tau_{T}}^{X_{\tau_{T}}^{\ast
,\pi
}(Z)}
\bigl(x,X_{\tau_{T}}^{\ast,\pi
}(Z)\bigr)\Gamma_{T}^{\ast}(Z)
\bigr].
\]
\end{theorem}


\begin{pf}
Using the definition of $X^{\ast,\pi}$ we
have for $g\in C_{c}^{\infty}(\mathbb{R}^d)$ (we recall that $\tau
_{T}\equiv\tau
_{J_{T}})$
\[
\mathbb{E} \bigl[ g\bigl(X_{T}^{*,\pi}\bigr)C_{T-\tau_{T}}(y)
\mid\tau_T, X_{\tau_{T}}^{*,\pi}=y \bigr] =\int
g(x)p_{T-\tau_{T}}^{\ast,y}(y,x)\,dx
\]
so that (\ref{Ad6}) says that $P_{T}^{\ast}(y,dx)=p_{T}^{\ast}(y,x)\,dx$
with
%
\begin{equation}
p_{T}^{\ast}(y,x)=e^{T}\mathbb{E} \bigl[
p_{T-\tau_{T}}^{\ast,X_{\tau
_{T}}^{\ast
,\pi}(y)}\bigl(X_{\tau_{T}}^{\ast,\pi}(y),x\bigr)
\Gamma_{T}^{\ast}(y) \bigr]. \label{Ad7}
\end{equation}
Notice that $p_{T}^{\ast}(y,x)=p_{T}(x,y)$ so the above equality says
that $%
P_{T}(x,dy)=p_{T}^{\ast}(y,x)\,dy$ with $p_{T}^{\ast}(y,x)$ given in the
previous formula. We conclude that the representation formula (\ref{Ad7})
proves that $P_{t}(x,dy)$ is absolutely continuous and the density is
represented by
\[
p_{T}(x,y)-p_{T}^{y}(x,y) =e^{T}
\mathbb{E} \bigl[ p_{T-\tau
_{T}}^{X_{\tau
_{T}}^{\ast,\pi}(y)}\bigl(x,X_{\tau_{T}}^{\ast,\pi}(y)
\bigr)\Gamma_{T}^{\ast
}(y)1_{\{J_{T}\geq1\}} \bigr].
\]
The representation for $P_{T}h$ can be obtained by integrating $\int
h(y)p_{T}(x,y)\,dy$ using (\ref{eq:adjden}).
\end{pf}

As before, we also have that the following generalized formulas with a
general Poisson process with parameter $\lambda$ are valid:
\[
p_{T}(x,y)-p_{T}^{y}(x,y)=e^{\lambda T}
\mathbb{E} \bigl[ \lambda ^{-J_{T}^{\lambda
}}p _{T-\tau_{J_{T}^{\lambda}}}^{X_{\tau_{J_{T}^{\lambda}}}^{\ast
,\pi}(y)}
\bigl(x,X_{\tau_{J_{T}^{\lambda}}}^{\ast
,\pi}(y)\bigr)\Gamma_{T}^{\ast}(y)1_{\{J_{T}^{\lambda}\geq1\}}
\bigr].
\]

We discuss now the regularity of $p_{t}(x,y)$.

\begin{theorem}
\label{b2}Suppose Hypotheses \ref{A0ast}, \ref{A1ast} and \ref{A2ast}.
\begin{longlist}[(ii)]
\item[(i)]
Furthermore assume that $(t,x,y)\rightarrow(p^y_t(x,y),\hat a_t(x,y))$
is continuous in $(0,T]\times\mathbb{R}^d\times\mathbb{R}^d$. Then
$(t,x,y)\rightarrow p_{t}(x,y)$ is continuous on $(0,\infty)\times
\mathbb{R}^d\times\mathbb{R}^d$. Moreover,
\[
p_{T}(x,y)=e^{T}\mathbb{E} \bigl[ p_{T-\tau_{_{T}}}^{X_{\tau
_{T}}^{\ast,\pi
}(y)}
\bigl(x,X_{\tau_{T}}^{\ast,\pi}(y)\bigr)\hat{\Gamma}_{T}(y)
\bigr]. \label{Ad8}
\]
\item[(ii)]
Furthermore, we assume that $x\rightarrow p_t^y(x,y)$ is differentiable
for all $(t,y)\in\mathbb{R}_+\times\mathbb{R}^d$ and that the
Hypothesis \ref{(h2a')} is satisfied for $\nabla_xp_t^y(x,y)$ instead
of $G$.
Then
the function $%
x\rightarrow p_{t}(x,z)$ is one time differentiable. Moreover,
\[
\nabla_{x}p_{T}(x,y)=\mathbb{E} \bigl[
\nabla_{x}p _{T-\tau
_{_{T}}}^{X_{\tau
_{T}}^{\ast,\pi}(y)}\bigl(x,X_{\tau
_{T}}^{\ast,\pi}(y)
\bigr)\hat{\Gamma}_{T}(y) \bigr]. \label{H6}
\]
\end{longlist}
\end{theorem}

\section{Examples: Applications to stochastic differential
equations}\label%
{sec:contdiff}

In this section, we will consider the first natural example for our previous
theoretical developments, that is, the case of multidimensional diffusion
processes. The forward method will need smooth coefficients and the
backward method will
require H\"older continuous coefficients.


\subsection{Example~1: The forward method for continuous SDE's with smooth
coefficients}
\label{sec:6.1}
We consider the following $d$-dimensional SDE:
%
\begin{equation}
\label{eq:SDE} X_{t}=x+\sum_{j=1}^{m}
\int_{0}^{t}\sigma _{j}(X_{s})\,dW_{s}^{j}+
\int_{0}^{t}b(X_{s})\,ds.
\end{equation}
Here, $\sigma_{j}$, $b\dvtx \mathbb{R}^d\rightarrow\mathbb{R}^d$, $\sigma
_{j}\in
C_{b}^{2}(\mathbb{R}^d;\mathbb{R}^d)$ is uniformly elliptic (i.e.,
$0<\underline{a}I\leq a\leq\overline{a}%
I $ for $\underline{a}$, $\overline{a}\in\mathbb{R}$ with $a=\sigma
\sigma^{\ast}$), $b\in C_{b}^{2}(\mathbb{R}^d;\mathbb{R}^d)$ and $W$
is a $m$-dimensional Wiener process. Under these conditions, there
exists a unique pathwise solution to the above equation. Then we define
the semigroup $%
P_{t}f(x)=\mathbb{E}[f(X_{t})]$ which has infinitesimal generator given
by $Lf(x)=\frac{1}{2}\sum_{i,j}a^{i,j}(x)\partial
_{i,j}^{2}f(x)+\sum_{i}b^{i}(x)\partial_{i}f(x)$ for $f\in
C_{c}^{\infty
}(\mathbb{R}^d)$ and $a^{i,j}(x)=\sum_k\sigma^i_k(x)\sigma^j_k(x)$.
Clearly, $P_tf(x)$ is jointly measurable and bounded and, therefore,
Hypothesis \ref{(H1)} is satisfied. We will consider the following approximation
process:
\[
X_{t}^{z}(x)=x+\sum_{j=1}^{m}
\sigma_{j}(z)W_{t}^{j}+b(z)t,
\]
which defines the semigroup
%
\begin{equation}
\label{eq:Pzdiff} P_{t}^{z}f(x)=\mathbb{E}\bigl[f
\bigl(X^z_t(x)\bigr)\bigr]=\int f(y)q_{ta(z)}
\bigl(y-x-b(z)t\bigr)\,dy,
\end{equation}
for $f\in C_c^\infty(\mathbb{R}^d)$, with jointly continuously
differentiable probability kernel $p_{t}^{z}(x,y)=q_{ta(z)}(y-x-b(z)t)$.
Furthermore, its associated infinitesimal operator [for $f\in
C_{c}^{2}(\mathbb{R}^d)
$] is given by
\[
L^{x}f(y)=\frac{1}{2}\sum_{i,j}a^{i,j}(x)
\partial _{i,j}^{2}f(y)+\sum_{i}b^{i}(x)
\partial_{i}f(y). \label{def:LXdif}
\]
Therefore, Hypotheses \ref{(H2)} and \ref{(H3)} are clearly satisfied.
Hypothesis \ref{(A0)} is clearly satisfied as $a^{i,j}, b^i\in
C_b^2(\mathbb{R}^d)$ for $i, j\in\{1,\ldots,d\}$. Now we proceed with the
verification of Hypothesis \ref{(A3)}.
Using integration by parts, we have for $f\in C_c^\infty(\mathbb{R}^d)$
%
\begin{eqnarray}\label{def:Sdiff}
S_{t}f(x) &=&\int\bigl(L^{y}-L^{x}
\bigr)f(y)P_{t}^{x}(x,dy) \nonumber
\\
&=&\frac{1}{2}\sum_{i,j}\int
\bigl(a^{i,j}(y)-a^{i,j}(x)\bigr)q_{ta(x)}
\bigl(y-x-b(x)t\bigr)\partial_{i,j}^{2}f(y)\,dy
\nonumber
\\
&&{}+\sum_{i}\int\bigl(b^{i}(y)-b^{i}(x)
\bigr)q_{ta(x)}\bigl(y-x-b(x)t\bigr)\partial_{i}f(y)\,dy
\\
&=&\int dy\,f(y) \biggl( \frac{1}{2}\sum_{i,j}
\partial _{i,j}^{2} \bigl(\bigl(a^{i,j}(y)-a^{i,j}(x)
\bigr)q_{ta(x)}\bigl(y-x-b(x)t\bigr) \bigr) \biggr)
\nonumber
\\
&&{}-\sum_{i}\int dy\,f(y)\partial_{i}
\bigl( \bigl(b^{i}(y)-b^{i}(x)\bigr)q_{ta(x)}
\bigl(y-x-b(x)t\bigr) \bigr).
\nonumber
\end{eqnarray}
In view of (\ref{G3}), we have
\begin{eqnarray*}
\partial_{i,j}^{2} \bigl( \bigl(a^{i,j}(y)-a^{i,j}(x)
\bigr)q_{ta(x)}\bigl(y-x-b(x)t\bigr) \bigr) &=&\theta_{t}^{i,j}(x,y)q_{ta(x)}
\bigl(y-x-b(x)t\bigr),
\\
\partial_{i} \bigl( \bigl(b^{i}(y)-b^{i}(x)
\bigr)q_{ta(x)}\bigl(y-x-b(x)t\bigr) \bigr) &=&\rho _{t}^{i}(x,y)q_{ta(x)}
\bigl(y-x-b(x)t\bigr),
\end{eqnarray*}
where we define for the Hermite polynomials $H$ (see Section~\ref{sec:Gaussian})
%
\begin{eqnarray}
 \theta_{t}^{i,j}(x,y) &=&\partial_{i,j}^{2}a^{i,j}(y)+
\partial _{j}a^{i,j}(y)h_{t}^{i}(x,y)+
\partial _{i}a^{i,j}(y)h_{t}^{j}(x,y)\nonumber\\
&&{}+
\bigl(a^{i,j}(y)-a^{i,j}(x)\bigr)h_{t}^{i,j}(x,y),
\nonumber
\\
\rho_{t}^{i}(x,y) &=&\partial_{i}b^{i}(y)+
\bigl(b^{i}(y)-b^{i}(x)\bigr)h_{t}^{i}(x,y),
\nonumber
\\
h_{t}^{i}(x,y)&=&H_{ta(x)}^{i}
\bigl(y-x-b(x)t\bigr),
\\
h_{t}^{i,j}(x,y)&=&H_{ta(x)}^{i,j}
\bigl(y-x-b(x)t\bigr).\label{D12}
\end{eqnarray}
So we obtain
%
\begin{eqnarray}\label{eq:thetadiff}
S_{t}f(x)&=&\int dy\,f(y)q_{ta(x)}\bigl(y-x-b(x)t\bigr)
\theta_{t}(x,y)
\nonumber
\\[-8pt]
\\[-8pt]
\nonumber
&=&\int f(y)\theta _{t}(x,y)P_{t}^{x}(x,dy).
\end{eqnarray}
Therefore, we have that
\[
\theta_{t}(x,y)=\frac{1}{2}\sum_{i,j}
\theta_{t}^{i,j}(x,y)-\sum_{i}
\rho _{t}^{i}(x,y).
 \label{D1'}
\]

Now, we verify Hypotheses \ref{(A3)} and \ref{A3p}. We have verified
the first part of Hypothesis \ref{(A3)} by the definition of $\theta$
in (\ref{eq:thetadiff}). In order to verify the rest of the conditions
in Hypothesis \ref{(A3)}, we see that by \eqref{G3} and (\ref{G4}) with
$\alpha=1$
\[
p_t^x(x,y)\bigl\llvert \theta_{t}(x,y)\bigr
\rrvert \leq C \bigl( \llVert a\rrVert _{2,\infty}+\llVert b\rrVert
_{1,\infty} \bigr) \frac
{1}{t^{{1}/2}%
}q_{ct\overline{a}}(y-x)
\]
for a constant $C>1$ and $c\in(0,1)$, and consequently all the
conditions in Hypothesis \ref{(A3)} are satisfied with $\rho=\frac{1}2+\rho_0$ and $%
\gamma_{t}(x,y)=\beta_t(x,y)=t^{\rho_0} q_{ct\overline{a}}(y-x)$.
Here, $\rho_0\in(\frac{\zeta-1}{2},\frac{1}2)$.

Similarly, Hypothesis \ref{A3p} is satisfied under the $\xi(x)=C$ for
\eqref{H2} by using that $1_{\{\sum_{i=1}^{n-1}\llvert  y_{i}\rrvert \leq
R\}}\le1$. For \eqref{h2a}, one uses that
\[
1_{\{\sum_{i=1}^{n-1}\llvert  y_{i}\rrvert >
R\}}\le\sum_{i=1}^{n-1}\sum_{j=1}^d1_{\{\llvert  y^j_{i}\rrvert >
{R}/{(n\sqrt{d})}\}}.
\]
Next, one performs the change of variables
$y_1=x_1$, $y_i-y_{i+1}=x_{i+1}$ for $i=1,\ldots,n-2$ in the integral of
\eqref{h2a} and use the inequality $1_{\{\llvert  y^j_{i}\rrvert >
{R}/{(n\sqrt{d})}\}}\le\frac{n^2d|y^j_i|^2}{R^2}$ to obtain the
following bound:
%
\begin{eqnarray}\quad
\label{eq:rost}&& \sum_{i=1}^{n-1}\sum
_{j=1}^d\frac{n^2d}{R^2}\sup
_{(y_0,y_n)\in
B(z_0,r)\times B(z_n,r)}\int dx_{1}\cdots\int dx_{n-1}\Biggl
\llvert \sum_{k=1}^ix^j_k
\Biggr\rrvert ^2 {q}_{c\delta_0\overline{a}}(x_1-y_0)
\nonumber
\\[-8pt]
\\[-8pt]
\nonumber
&&\qquad{}\times\prod_{i=1}^{n-2}{q}_{c\delta_i\overline
{a}}(x_{i+1}){q}_{c\delta_{n-1}\overline{a}}
\Biggl(y_n-\sum_{i=1}^{n-1}x_i
\Biggr).
\end{eqnarray}
Without loss of generality, using a further change of variables
$z_1=x_1-y_0$, we may consider the case where $y_0=0$. Next, we use the
inequality
$\llvert \sum_{k=1}^ix^j_k\rrvert ^2\le n\sum_{k=1}^i|x^j_k|^2$. Then
one rewrites the integral in a probabilistic way using Gaussian random
variables. This becomes $\mathbb
{E}[|Z^j_k|^2/Z_0+\cdots+Z_{n-1}=y_n]p_{Z_0+\cdots+Z_{n-1}}(y_n)$ where $Z_i$
is a $d$-dimensional Gaussian random vector with mean 0 and covariance
matrix $c\delta_i\overline{a}I$. The conditional variance can be computed
explicitly and the density can be bounded by its maximum value (i.e.,
$y_n=0$). Finally, we obtain that \eqref{eq:rost} is bounded by
$C\frac{n^4d^2}{R^2\sqrt{\delta}} (|y_n|^2+\delta)
$ with $\delta=\sum_{i=0}^{n-1}\delta_i$.
Therefore, condition \eqref{h2a} will be satisfied taking
$R_\varepsilon=\varepsilon^{-{1}/2}$ and $\xi(\delta)=\delta^{-{1}/2}+\delta^{{1}/2}$ and the upper bound in \eqref{h2a} becomes
$C2^n\xi(\delta)$.
We leave the details of the calculation for the reader.
Therefore, the existence of the density follows.

In order to obtain further regularity, we need to verify the uniform
integrability condition for $\zeta\in(1,\rho^{-1})$. In this case, we
first note that due to \eqref{G4}, $|\nabla_yp_t^x(x,y)|\le\frac{C}{t^{{1}/2}}q_{ct\overline{a}}(y-x)$. Therefore, we may choose any
$\rho\in(\frac{1}2,\frac{2}3)$ and let $\zeta=\frac{1}{2(1-\rho)}> 1$.
Finally, we define
$\gamma_t(x,y)=t^{({1}/2)(1-\zeta^{-1})}q_{ct\overline{a}}(y-x)$ and
$\xi(x)=C$ in order to obtain \eqref{H2}. One also obtains \eqref{h2a}
as in the proof of continuity. Therefore, the hypotheses in Theorem~\ref
{f4} are satisfied.

Now, we give the description of the stochastic representation. Given a
Poisson process with parameter $\lambda=1$ and jump times
$\{\tau_i, i=0,\ldots\}$. Given that $J_T=n$ and $t_i:=\tau_i\wedge T$ we
define the
process $ ( X_{t_{i}}^{\pi} ) _{i=0,\ldots,n+1}$ for $\pi
=\{t_{i};i=0,\ldots,n+1\}$, with $0=t_{0}<t_{1}<\cdots<t_{n}\leq t_{n+1}=T$ is
then defined as compositions of $X^{z}(x)$ as follows:
\[
X_{t_{k+1}}^{\pi}= X_{t_{k+1}-t_{k}}^{z}(x)
\vert _{z=x=X_{t_{k}}^{\pi}},
\]
for $k=0,\ldots,n$. Here $X_{0}^{\pi}=x$ and the noise used for $%
X_{t_{k+1}-t_{k}}^{z}(x)$ is independent of $X_{t_{j}}^{\pi}$ for all $
j=0,\ldots,k$ and of the Poisson process $J$.

\begin{theorem}
Suppose that $a\in C_{b}^{2}(\mathbb{R}^d;\mathbb{R}^d\times\mathbb
{R}^d),b\in
C_{b}^{2}(\mathbb{R}^d;\mathbb{R}^d)$ and $\overline{a}\geq a\geq
\underline{a}$. Define
\[
\Gamma_{T}(x)= \cases{ %
\displaystyle\prod
_{j=0}^{J_{T}-1}\theta_{\tau_{j+1}-\tau_{j}}
\bigl(X_{\tau_{j}}^{\pi
},X_{\tau_{j+1}}^{\pi}\bigr), & \quad $\mbox{if }J_{T}\geq1,$
\vspace*{2pt}\cr
1, &\quad  $\mbox{if }J_{T}=0.$}
\]
Then for any $f\in C_c^\infty(\mathbb{R}^d)$ we have
\[
P_{T}f(x)=e^{T}\mathbb{E} \bigl[f\bigl(X_T^\pi
\bigr)\Gamma _{T}(x) \bigr]
\]
and, therefore,
\[
p_{T}(x,y)=e^{T}\mathbb{E} \bigl[ p _{T-\tau_{T}}^{X_{\tau_{T}}^{\pi
}}
\bigl(X_{\tau_{T}}^{\pi},y\bigr)\Gamma _{T}(x) \bigr],
\label{D4}
\]
where $ ( X_{t}^{\pi} ) _{t\in\pi}$ is the Euler scheme
with $%
X_{0}^{\pi}=x$ and random partition $\pi=\{\tau_{i};i=0,\ldots,\tau
_{J_{T}}\}\cup\{T\}$ where $0=\tau_{0}<\cdots<\tau_{J_{T}}\leq T$ where the
random times $\{\tau_{i}\}_i$ are the associated jump times of the
simple Poisson
process $J$, independent of $X^\pi$ with $\mathbb{E}[J_{T}]=T$.
Moreover, $(t,x,y)\rightarrow p_{t}(x,y)$ is
continuous on $(0,\infty)\times\mathbb{R}^d\times\mathbb{R}^d$ and
for every $t>0$ the
function $(x,y)\rightarrow p_{t}(x,y)$ is continuously differentiable. We
also have
\[
\partial_{y^{i}}p_{T}(x,y)=e^{T}\mathbb{E} \Biggl[
h_{T-\tau
_{T}}^{i}\bigl(X_{\tau
_{T}}^{\pi},y\bigr)p
_{T-\tau_{T}}^{X_{\tau_{T}}^{\pi}}\bigl(X_{\tau
_{T}}^{\pi
},y\bigr)
\prod_{j=0}^{J_{T}-1}\theta_{\tau_{j+1}-\tau_{j}}
\bigl(X_{\tau
_{j}}^{\pi
},X_{\tau_{j+1}}^{\pi}\bigr)
\Biggr], \label{D5}
\]
where $h^{i}$ is defined in (\ref{D12}).
\end{theorem}

\begin{pf}
As a consequence of Theorems \ref{th:forwardrep}, \ref{f3} and \ref
{f4}, we obtain most of the
mentioned results. The fact that $y\rightarrow
p_{t}(x,y) $ is continuously differentiable will follow from the
backward method concerning the adjoint semigroup that we present in the
following
section (since $a$ is differentiable it is also H\"{o}lder continuous
so the
hypotheses in the next section are verified).
\end{pf}

\subsection{Example~2: The backward method for continuous SDEs with H\"{o}lder
continuous coefficients}

\label{sec:6.2}
In this section, we will assume the same conditions as in the previous
section except the regularity hypothesis on $a$ and $b$. We will assume
that $a$ is a H\"{o}lder continuous function of
order $\alpha\in(0,1)$ and $b$ is a bounded measurable function.
We suppose the existence of a unique weak solution to \eqref{eq:SDE}.
For further references on this matter, see \cite{SV}. The
approximating semigroup is the same as in the previous section and is given
by (\ref{eq:Pzdiff}). Therefore we have, as before,
\begin{eqnarray*}
p_{t}^{z}(x,y) &=&q_{ta(z)}\bigl(y-x-b(z)t\bigr),
\\
\phi_{t}^{z}(x) &=&q_{ta(z)}\bigl(z-x-b(z)t\bigr).
\end{eqnarray*}
In this case, note that for fixed $z\in\mathbb{R}^d$, $\phi^{z}$ is a
smooth density function
and therefore $C_{t}(x)=1$. Furthermore, as in the previous section,
Hypotheses \ref{(H1)}, \ref{(H2)} and~\ref{(H3)} are satisfied.
Similarly, Hypothesis \ref{A0ast} can be easily verified. We will now
check Hypothesis \ref{A1ast}. We define
\begin{eqnarray*}
\hat{\theta}_{t}(x,z) &=&\frac{1}{2}\sum
_{i,j}\bigl(a^{i,j}(x)-a^{i,j}(z)\bigr)
\hat{h}%
_{t}^{i,j}(x,z)-\sum
_{i}\bigl(b^{i}(x)-b^{i}(z)\bigr)
\hat{h}_{t}^{i}(x,z),
\\
\hat{h}_{t}^{i}(x,z) &=&H_{ta(z)}^{i}
\bigl(z-x-b(z)t\bigr),
\\
\hat{h}_{t}^{i,j}(x,z) &=&H_{ta(z)}^{i,j}
\bigl(z-x-b(z)t\bigr)
\end{eqnarray*}
so that, by (\ref{G3}),
\begin{eqnarray*}
\bigl(L^{x}-L^{z}\bigr)\phi_{t}^{z}(x)
&=&\frac{1}{2}%
\sum_{i,j}
\bigl(a^{i,j}(x)-a^{i,j}(z)\bigr)\partial _{i,j}^{2}q_{ta(z)}
\bigl(z-x-b(z)t\bigr)\\
&&{} -\sum_{i}\bigl(b^{i}(x)-b^{i}(z)
\bigr)\partial _{i}q_{ta(z)}\bigl(z-x-b(z)t\bigr)
\\
&=&\hat{\theta}_{t}(x,z)q_{ta(z)}\bigl(z-x-b(z)t\bigr).
\end{eqnarray*}
Using (\ref{G4}) and the H\"{o}lder continuity of $a^{i,j}$, we obtain
\begin{eqnarray*}
&&\bigl\llvert \bigl(a^{i,j}(x)-a^{i,j}(z)\bigr)
\partial_{i,j}^{2}\phi_{t}^{z}(x)\bigr
\rrvert\\
&&\qquad\leq C\llvert x-z\rrvert ^{\alpha}\bigl\llvert
\partial_{i,j}^{2}\phi _{t}^{z}(x)\bigr
\rrvert\\
&&\qquad \leq C\bigl(\bigl\llvert z-x-b(z)t\bigr\rrvert ^{\alpha}+\llVert b
\rrVert _{\infty}^{\alpha}t^{\alpha}\bigr)\bigl\llvert \partial
_{i,j}^{2}q_{ta(z)}\bigl(z-x-b(z)t\bigr)\bigr\rrvert
\\
&&\qquad\leq Ct^{-(1-\alpha/2)}q_{\overline{a}t}\bigl(z-x-b(z)t\bigr).
\end{eqnarray*}
And using (\ref{G4})(ii) with $\alpha=0$, we obtain
\[
\bigl\llvert \bigl(b^{i}(x)-b^{i}(z)\bigr)
\partial_{i}\phi_{t}^{z}(x)\bigr\rrvert \leq
\frac{2}{%
t^{{1}/2}}\llVert b\rrVert _{\infty}q_{t\overline{a}}\bigl(z-x-b(z)t
\bigr).
\]
Finally, we have
\[
\bigl\llvert \hat{\theta}_{t}(x,z)\bigr\rrvert \leq\frac{C}{t^{1-\alpha/
2}}
\bigl(1+\llVert b\rrVert _{\infty}\bigr)q_{t\overline{a}}\bigl(z-x-b(z)t
\bigr).
\]
We also have $\phi_{t}^{z}(x) \leq Cq_{t\overline{a}%
}(z-x-b(z)t)$ so we obtain
\[
\phi_{t}^{z}(x)\bigl\llvert \hat{\theta}_{t}(x,z)
\bigr\rrvert \leq\frac{C}{%
t^{1-\alpha/2}}\bigl(1+\llVert b\rrVert _{\infty}
\bigr)q_{2t\overline{a}%
}\bigl(z-x-b(z)t\bigr).
\]
We conclude that Hypothesis \ref{A1ast} is verified. The verification
of Hypothesis \ref{A2ast} is done like in the previous section using
$\rho\in(\frac{2-\alpha}2,\frac{3-\alpha}3)$ and $\zeta=(3-\alpha-2\rho
)^{-1}\in(1,\rho^{-1})$. Therefore, we have the following result.

\begin{proposition}
Suppose that $a$ is H\"{o}lder continuous of order $\alpha\in(0,1)$, $
\overline{a}\geq a\geq\underline{a} $ and $b$ is measurable and
bounded. Then
\[
p_{T}(x,y)=e^{T}\mathbb{E} \Biggl[ p _{T-\tau_{T}}^{X_{\tau_{T}}^{\ast
,\pi
}(y)}
\bigl(x,X_{\tau_{T}}^{\ast,\pi
}(y)\bigr)\prod
_{j=0}^{J_{T}-1}\hat{\theta}_{\tau_{j+1}-\tau_{j}}
\bigl(X_{\tau
_{j+1}}^{\ast,\pi}(y),X_{\tau_{j}}^{\ast,\pi}(y)
\bigr) \Biggr], \label{B2}
\]
where $X^{\ast,\pi}(y)$ is the Euler scheme with $X_{0}^{\ast,\pi
}=y$ and drift coefficient $-b$.
Moreover, $(t,x,y)\rightarrow p_{t}(x,y)$ is continuous on $(0,\infty
)\times
\mathbb{R}^d\times\mathbb{R}^d$ and for every $(t,y)\in(0,\infty
)\times\mathbb{R}^d$ the
function $x\rightarrow p_{t}(x,y)$ is continuously differentiable.
Moreover,
\begin{eqnarray*}
&&\partial_{x^{i}}p_{T}(x,y)\\
&&\qquad=-e^{T}\mathbb{E}
\Biggl[ \hat{h}_{T-\tau
_{T}}^{i}\bigl(x,X_{\tau_{T}}^{\ast,\pi}(y)
\bigr)p _{T-\tau_{T}}^{X_{\tau
_{T}}^{\ast,\pi}(y)}\bigl(x,X_{\tau_{T}}^{\ast,\pi
}(y)
\bigr)\\
&&\hspace*{50pt}\qquad\quad{}\times\prod_{j=0}^{J_{T}-1}\hat{
\theta}_{\tau_{j+1}-\tau
_{j}}\bigl(X_{\tau_{j+1}}^{\ast,\pi}(y),X_{\tau_{j}}^{\ast,\pi
}(y)
\bigr) \Biggr].
\end{eqnarray*}
\end{proposition}

\subsection{Example~3: One-dimensional L\'{e}vy driven SDE with H\"{o}lder
type coefficients}
\label{sec:jumps}

Although we may consider various other situations where the forward and the
backward method can be applied and to test their limits, we prefer to
concentrate in this section on the backward method for a one-dimensional
jump type SDEs driven by a L\'{e}vy process of a particular
type: we assume that the intensity measure of the L\'{e}vy process is a
mixture of Gaussian densities. This a quite general class as it can be
verified from Schoenberg's theorem; see \cite{Re}.

For this, let $N(dx,dc,ds)$ denote the Poisson random measure associated
with the compensator given by $q_{c}(x)\,dx\,\nu(dc)\,ds$ where $\nu$
denotes a
nonnegative measure on $\mathbb{R}_+:=(0,\infty)$ which satisfies the
following.

\begin{guess}
\label{(P1)}
%
$\nu(\mathbb{R}_+)=\infty$ and $C_{\nu}:=\int_{\mathbb{R}_+}c\nu
(dc)<\infty$.
%
\end{guess}

We refer the reader to \cite{IW} for notation and
detailed definitions on Poisson random measures.
Therefore, heuristically speaking, $x$ stands for the jump size which
arises from a Gaussian distribution with random variance obtained from the
measure $\nu$.

We define $\eta_{\nu}(u):=\nu(u,\infty)$ and we assume that there exists
some $s_{\ast}\geq0$ and $h,C_{\ast}>0$ such that we have the
following.%

\begin{guess}
\label{(P2)}
%
$\int_0^\infty e^{-u}\eta_{\nu}(\frac{u}s)\,du\geq C_{\ast}s^{h}\int_0^\infty e^{-u}\eta_{\nu}(u)\,du$
$\forall s\geq s_{\ast}$.
%
\end{guess}

For example, if $\nu(dc)=1_{(0,1]}(c)c^{-(1+\beta)}\,dc$ with $0<\beta<1$
then Hypothesis~\ref{(P1)} is satisfied and Hypothesis \ref{(P2)} is
satisfied with $h=\beta$.

$\widetilde{N}(dx,dc,ds)=N(dx,dc,ds)-q_{c}(x)\,dx\,\nu(dc)\,ds$ denotes the
compensated Poisson random measure. We also define the following auxiliary
processes and driving process $Z$:
\begin{eqnarray*}
V_{t} &=&\int_{0}^{t}\int
_{\mathbb{R}_+\times\mathbb{R}}cN(dx,dc,ds),
\\
Z_{t} &=&\int_{0}^{t}\int
_{\mathbb{R}_+\times\mathbb{R}}xN(dx,dc,ds),
\\
N_{\nu}(dx,ds) &=&\int_{\mathbb{R}_+}N(dx,dc,ds).
\end{eqnarray*}
With a slight variation of some classical proofs (see, e.g., Chapter~2
in \cite{Appl}) one can obtain the following generalization of the L\'
evy--Khinchine formula.

\begin{proposition}
\label{prop:chfL}
Assume Hypothesis \ref{(P1)}. Let $h\dvtx \mathbb{R}\times\mathbb
{R}_+\rightarrow\mathbb{R}$ be such that $\llvert \int_{\mathbb
{R}\times\mathbb{R}_+}(e^{i\theta h(x,c)}-1)q_c(x)\,dx\,d\nu(c)\rrvert <\infty$. Then the stochastic process $U_t(h):=\int_{0}^{t}\int_{\mathbb{R}_+\times\mathbb{R}}h(x,c)N(dx,dc,ds)$ has independent
increments with characteristic function given by
\[
\mathbb{E}\bigl[\exp\bigl(i\theta U_t(h)\bigr)\bigr]=\exp \biggl(t
\int_{\mathbb{R}\times
\mathbb{R}_+}\bigl(e^{i\theta h(x,c)}-1\bigr)q_c(x)\,dx\,d
\nu(c) \biggr)
\]
%
the density of $Z_t$ at $y$ can be written as
$
\mathbb{E}[q_{V_t}(y)]$.
\end{proposition}

\begin{pf}
The first part of the proof is classical, while in order to obtain the
representation for the density of $Z_t$, one takes $h(x,c)=x$ to obtain
the characteristic function associated with $Z_t$ under Hypothesis \ref
{(P1)}. On the other hand, one only needs to compute the characteristic
function associated with the density function $\mathbb{E}[q_{V_t}(y)]$
to finish the proof.
\end{pf}

Notice that due to Hypothesis \ref{(P1)} we have that
%
\begin{equation}
\label{eq:2mom} \mathbb{E}\bigl[Z_{t}^{2}\bigr]=t\int
_{\mathbb{R}\times\mathbb{R}_+}\llvert u\rrvert ^{2}q_{c}(u)\nu
(dc)\,du=t\int_{\mathbb{R}_+}c\nu(dc)<\infty.
\end{equation}
Therefore, $Z$ is a L\'{e}vy process of finite variance. $N_{\nu
}(dx,ds)$ is
a Poisson random measure with compensator $\mu_\nu(dx)\,ds:=\int_{\mathbb
{R}_+}q_{c}(x)\nu(dc)\,dx\,ds$
and we denote by $\widetilde{N}_{\nu}(dx,ds)$ the compensated Poisson
random measure. Then we consider the
solution of the following stochastic differential equation driven by
$Z$ and
its corresponding approximation obtained after freezing the jump
coefficient. That is,\vspace*{-6pt}
{\renewcommand{\theequation}{$E_{\nu}$}
\begin{equation}\label{enu}
 X_{t}^{\nu}(x) =x+\int
_{0}^{t}\int_{\mathbb
{R}}\sigma
\bigl(X_{s-}^{\nu}(x)\bigr)u\widetilde{N}_{\nu}(ds,du),\vspace*{-6pt}
\end{equation}}
{\renewcommand{\theequation}{$E_{\nu}^{z}$}
\begin{equation}\label{enuz}
 X_{t}^{\nu,z}(x)
=x+\int_{0}^{t}\int_{\mathbb
{R}}
\sigma(z)u%
\widetilde{N}_{\nu}(ds,du).
\end{equation}}
\hspace*{-4pt}We assume that $\sigma\dvtx\mathbb{R}\rightarrow\mathbb{R}$ verifies the
following conditions.%

\begin{guess}
\label{(P3)}
%
\textup{(i)} There exists $\underline{\sigma}, \overline{\sigma}%
>0$ such that $\underline{\sigma} \leq\sigma(x)\leq\overline{%
\sigma}$ for all $x\in\mathbb{R}$.

\textup{(ii)} There exists $\alpha\in(0,1]$ such that $\llvert  \sigma
(x)-\sigma(y)\rrvert  \leq C_{\alpha}\llvert
x-y\rrvert  ^{\alpha}$.
%
\end{guess}

If $\alpha=1$, then \eqref{enu} has a unique solution. Here, rather
than entering into the discussion of existence and uniqueness results for
other values of $\alpha\in(0,1] $, we refer the reader to a survey
article by
Bass and the references therein (see \cite{B}). Therefore, from now on,
we suppose that a
unique weak solution to \eqref{enu} exists so that $P_{t}^{\nu
}f(x)=\mathbb{E} [ f(X_{t}^{\nu
}(x)) ] $ is a semigroup with infinitesimal operator [note that
$\int u\mu_\nu(du)=0$]
\[
L^{\nu}f(x)=\int_{\mathbb{R}} \bigl(f\bigl(x+\sigma(x)u
\bigr)-f(x)\bigr)\mu_{\nu
}(du).
\]
Therefore, Hypothesis \ref{(H1)} is clearly satisfied.

Similarly, $X^{\nu,z}(x)$, defines a semigroup $P_{t}^{\nu
}f(x)=\mathbb{E} [
f(X_{t}^{\nu,z}(x)) ] $ with infinitesimal operator
%
\renewcommand{\theequation}{\arabic{section}.\arabic{equation}}
\setcounter{equation}{8}
\begin{equation}
L^{\nu,z}f(x)=\int_{\mathbb{R}} \bigl(f\bigl(x+\sigma(z)u
\bigr)-f(x)\bigr)\mu _{\nu}(du). \label{eq:genstar}
\end{equation}
Our aim is to give sufficient conditions in order that the law of $%
X_{t}^{\nu}(x)$ is absolutely continuous with respect to the Lebesgue
measure and to represent the density $p_{t}(x,y)$ using the backward method
as introduced in Section~\ref{sec:back}. In order to proceed with the
verification of Hypothesis \ref{(H2)}, we need to prove the following
auxiliary lemma.\vadjust{\goodbreak}

\begin{lemma}
Suppose that Hypotheses \ref{(P1)} and \ref{(P2)} holds for some $h>0$.
Then for every
$p>0$ there exists a constant $C$ such that for every $t>0$
%
\begin{equation}
\mathbb{E} \bigl[ V_{t}^{-p} \bigr] \leq
Ct^{-{p}/h}. \label{P2}
\end{equation}
\end{lemma}

\begin{pf} Recall that the Laplace transform of $V_{t}$ is given by
\[
\mathbb{E} \bigl[ e^{-aV_{t}} \bigr] =\exp \biggl( -t\int
_{\mathbb
{R}_+}\bigl(1-e^{-ac}\bigr)\nu (dc) \biggr).
\]
We use the change $s^{\prime}=sV_{t}$ and we obtain
\[
\int_{0}^{\infty}s^{p-1}e^{-sV_{t}}\,ds=c_{p}V_{t}^{-p}
\]
with $c_{p}=\int_{0}^{\infty}s^{p-1}e^{-s}\,ds$. It follows that
\[
c_{p}\mathbb{E} \bigl[ V_{t}^{-p} \bigr] =\int
_{0}^{\infty
}s^{p-1}\mathbb{E} \bigl[
e^{-sV_{t}}%
 \bigr] \,ds=\int_{0}^{\infty}s^{p-1}
\exp \biggl( -t\int_{\mathbb
{R}_+}\bigl(1-e^{-sc}\bigr)\nu
(dc) \biggr) \,ds.
\]
For $s>s_{\ast}$ we have using the integration by parts formula and the
change of variables $sc=u$,
\[
\int_{0}^{\infty}\bigl(1-e^{-sc}\bigr)
\nu(dc)=\int_{0}^{\infty}\,due^{-u}\eta
_{\nu
}\biggl(\frac{u}s\biggr)\geq C_{\ast}s^{h}
\int_{0}^{\infty}\,due^{-u}
\eta_{\nu
}(u)=:s^{h}\alpha_{\nu}
\]
with $\alpha_{\nu}\in\mathbb{R}_+$. Therefore, again by change of
variables, we
have that
\begin{eqnarray*}
\int_{s_{\ast}}^{\infty}s^{p-1}\exp \biggl( -t
\int_{\mathbb
{R}_+}\bigl(1-e^{-sc}\bigr)\nu (dc) \biggr) \,ds
&\leq&\int_{s_{\ast}}^{\infty}s^{p-1}e^{-ts^{h}\alpha
_{\nu
}}\,ds
\\
&\leq& t^{-{p}/h}C(\nu,p,h)
\end{eqnarray*}
with
\[
C(\nu,p,h)=h^{-1}\int_{0}^{\infty}u^{-(1-{p}/h)}e^{-u\alpha_{\nu
}}\,du<
\infty.
\]
Since $\int_{0}^{s_{\ast}}s^{p-1}\,ds=\frac{1}{p}s_{\ast}^{p}$ the
conclusion follows by taking $s_{\ast}=t^{-{1}/h}$.
\end{pf}

Now we can verify Hypothesis \ref{(H2)}. For this, we need to compute
as explicitly as possible the density $p_{t}^{z}(x,y)$
of the law of $X_{t}^{\nu,z}(x)$. In fact, the following is a
corollary of Proposition~\ref{prop:chfL} and the previous lemma which
is used together with Lemma~\ref{lem:estforq} in order to obtain the
needed uniform integrability properties.

\begin{corollary}
\label{lem:5.1}Suppose that Hypotheses \ref{(P1)} and \ref{(P3)} are
verified. Then the
law of $X_{t}^{\nu,z}(x)$ is absolutely continuous with respect to the
Lebesgue measure with strictly positive continuous density given by
\[
p_{t}^{z}(x,y)=\mathbb{E} \bigl[ q_{\sigma^{2}(z)V_{t}}(x-y)
\bigr]. \label{P1}
\]
Therefore, for each fixed $(t,z)\in(0,T]\times\mathbb{R}$, we have
that $p_t^z\in C^2_b(\mathbb{R}\times\mathbb{R})$ and $p_t^z(x,y)$ is
jointly continuous in $(t,z,x,y)$.
\end{corollary}

Note that due to the above result Hypothesis \ref{(H2)} is satisfied
and $\phi_{t}^{y}(x)=\mathbb{E} [ q_{\sigma
^{2}(z)V_{t}}(x-y) ]$.
Furthermore, as it is usually the case Hypotheses \ref{(H3)} and \ref
{A0ast}(0) are trivially satisfied.
For Hypothesis~\ref{A0ast}(i), one only needs to apply Corollary~\ref
{lem:5.1}. Hypothesis~\ref{A0ast}(ii) follows from the joint
continuity of $p_t^z(x,y)$ and Hypothesis~\ref{A0ast}(iii) follows
from the regularity of $p_t^z(x,y)$ as stated in the above Corollary~\ref{lem:5.1} and \eqref{eq:2mom}.

We are now ready to proceed and verify Hypotheses \ref{A1ast} and \ref
{A2ast}. We have by~(\ref{eq:genstar}), $\int uq_{c}(u)\,du=0$
and properties of convolution that
\begin{eqnarray*}
\bigl(L^{\nu,x}-L^{\nu,z}\bigr)\phi_{t}^{z}(x)
&=&\int_{\mathbb{R}_+\times
\mathbb{R}}\bigl(\phi _{t}^{z}
\bigl(x+\sigma(x)u\bigr)-\phi_{t}^{z}\bigl(x+\sigma(z)u
\bigr)\bigr)q_{c}(u)\nu(dc)\,du
\\
&=&\int_{\mathbb{R}_+\times\mathbb{R}} \bigl( \mathbb{E} \bigl[ q_{\sigma^{2}(z)V_{t}}
\bigl(x-z+\sigma (x)u\bigr) \bigr] \\
&&\hspace*{33pt}{}-\mathbb{E} \bigl[ q_{\sigma^{2}(z)V_{t}}\bigl(x-z+
\sigma (z)u\bigr) \bigr] \bigr) q_{c}(u)\nu(dc)\,du
\\
&=&\int_{\mathbb{R}}\mathbb{E} \bigl[ q_{\sigma^{2}(x)c+\sigma
^{2}(z)V_{t}}(x-z)-q_{\sigma
^{2}(z)c+\sigma^{2}(z)V_{t}}(x-z)
\bigr] \nu(dc).
\end{eqnarray*}
In particular, Hypothesis \ref{A1ast} holds with
%
\begin{eqnarray}\label{eq:thetajumps}
\hat{\theta}_{t}(x,y)&=&\frac{1}{\mathbb{E} [ q_{\sigma
^{2}(y)V_{t}}(x-y) ] }%
\nonumber
\\[-8pt]
\\[-8pt]
\nonumber
 &&{}\times\biggl\{ \int
_{\mathbb{R}_+}\mathbb{E} \bigl[ q_{\sigma^{2}(x)c+\sigma
^{2}(y)V_{t}}(x-y)-q_{\sigma^{2}(y)c+\sigma^{2}(y)V_{t}}(x-y)
\bigr] \nu (dc) \biggr\}.
\end{eqnarray}

\begin{theorem}
\label{th:74}
Suppose that Hypotheses \ref{(P1)}, \ref{(P2)} and \ref{(P3)} hold with
$h>$ $1-\frac\alpha2$.
Then the law of $X_{T}^{\nu}(x)$ is absolutely continuous with respect to
the Lebesgue measure and its density $p_{T}(x,y)$ satisfies
\[
p_{T}(x,y)=e^{T}\mathbb{E} \Biggl[ p_{T-\tau_{_{T}}}^{X_{\tau
_{T}}^{\ast,\pi
}(y)}
\bigl(x,X_{\tau_{T}}^{\ast,\pi}(y)\bigr)\prod
_{j=0}^{J_{T}-1}\hat{\theta }_{\tau
_{j+1}-\tau_{j}}
\bigl(X_{\tau_{j+1}}^{\ast,\pi}(y),X_{\tau_{j}}^{\ast
,\pi}(y)
\bigr)%
 \Biggr], \label{eq:pLevy}
\]
where $X_{t}^{\ast,\pi}(y)$ is the Euler scheme given in the backward
method starting at $X_0^{\ast, \pi}(y)=y$.
Moreover, $(t,x,y)\rightarrow p_{t}(x,y)$ is continuous on $(0,\infty
)\times
\mathbb{R}\times\mathbb{R}$ and for every $(t,y)\in(0,\infty)\times
\mathbb{R}$ the function $%
x\rightarrow p_{t}(x,y)$ is differentiable and
\[
\partial_{x}p_{T}(x,y)=e^{T}\mathbb{E} \Biggl[
\partial_{x}p_{T-\tau
_{_{T}}}^{X_{\tau_{T}}^{\ast,\pi}(y)}\bigl(x,X_{\tau_{T}}^{\ast,\pi
}(y)
\bigr)\prod_{j=0}^{J_{T}-1}\hat{
\theta}_{\tau_{j+1}-\tau_{j}}\bigl(X_{\tau
_{j+1}}^{\ast,\pi}(y),X_{\tau_{j}}^{\ast,\pi}(y)
\bigr) \Biggr].
\]
\end{theorem}

\begin{pf} We have already verified Hypotheses \ref{A0ast},
\ref{A1ast} and the differentiability of $p^z_t$. It remains to verify
the hypotheses in Lemma~\ref{lemma:1.10} and Theorem~\ref{th:1.13}. For
this, we have to estimate
\[
\bigl\llvert \hat{\theta}_{t}(x,y)\bigr\rrvert \phi_{t}^{y}(x)
\leq\int_{\mathbb
{R}_+}\mathbb{E} \bigl[ \bigl\llvert
q_{\sigma^{2}(x)c+\sigma^{2}(y)V_{t}}(x-y)-q_{\sigma
^{2}(y)c+\sigma^{2}(y)V_{t}}(x-y)\bigr\rrvert \bigr] \nu(dc).
\]
Let us denote $a=x-y$ and
\[
s^{\prime}=\sigma^{2}(y)c+\sigma^{2}(y)V_{t},\qquad
s^{\prime\prime
}=\sigma^{2}(x)c+\sigma^{2}(y)V_{t}.
\]
We assume that $s'\le s''$ (the other case is similar) and note the
inequality $($with $a,b,b^{\prime},c>0)$
\[
b\geq b^{\prime} \quad\Rightarrow\quad \frac{a+cb}{a+cb^{\prime
}}\leq \frac{b}{b^{\prime}}.
\]
From this inequality, we obtain
%
\begin{eqnarray}\label{eq:B}
\frac{s^{\prime\prime}}{s^{\prime}} &\leq&\frac{\overline{\sigma
}^{2}}{%
\underline{\sigma}^{2}} \quad\mbox{and}
\nonumber
\\[-8pt]
\\[-8pt]
\nonumber
\frac{\llvert  s^{\prime\prime}-s^{\prime}\rrvert  }{s^{\prime}} &\leq &%
\frac{c\llvert  \sigma^{2}(x)-\sigma^{2}(y)\rrvert  }{\sigma
^{2}(x)c+\sigma
^{2}(y)V_{t}}\leq\frac{cC_{\alpha}\llvert  a\rrvert  ^{\alpha
}}{\underline{%
\sigma}^{2}(c+V_{t})},
\end{eqnarray}
where $C_{\alpha}$ is the H\"{o}lder constant of $\sigma^{2}$. Finally,
from Lemma~\ref{lem:ineqgauss} and (\ref{eq:B}) this gives
\begin{eqnarray*}
\bigl\llvert q_{s^{\prime\prime}}(a)-q_{s^{\prime}}(a)\bigr\rrvert &\leq
&CC_{\alpha}\frac{\overline{\sigma}^{2}}{\underline{\sigma}^{2}}
\frac{c\llvert  a\rrvert  ^{\alpha}}{\underline{\sigma}^{2}(c+V_{t})}
q_{\overline{\sigma}^{2}(c+V_{t})}(a) \\
&\leq&
CC_{\alpha}\frac{\overline
{\sigma}^{2+\alpha}}{\underline{\sigma}%
^{4}} \frac{c}{(c+V_{t})^{1-\alpha/2}} q_{({\overline{%
\sigma}^{2}}/{2})(c+V_{t})}(a).
\end{eqnarray*}
Returning to our main proof, we obtain (with $C=CC_{\alpha}\overline
{\sigma}^{2+\alpha}%
\underline{\sigma}^{-4})$
%
\begin{eqnarray}\label{eq:TP}
&&\bigl\llvert \hat{\theta}_{t}(x,y)\bigr\rrvert \phi_{t}^{y}(x)
\nonumber
\\[-8pt]
\\[-8pt]
\nonumber
&&\qquad\leq C\int_{\mathbb
{R}_+}\mathbb{E} \bigl[ (c+V_{t})^{-(1-\alpha/2)}
q_{({\overline{\sigma}^{2}}/{2})(c+V_{t})}(x-y)%
 \bigr] c\nu(dc).
\end{eqnarray}

A first
step is to obtain estimates for the right-hand side of the above
inequality, so as to be able to define ${\gamma}$. For this, we define
%
\begin{eqnarray}\label{P6}
g_{t}(x,y) &=&\int_{\mathbb{R}_+}\mathbb{E} \bigl[
V_{t}^{-(1-
\alpha/2)}q_{({\overline{%
\sigma}^{2}}/{2})(c+V_{t})}(x-y) \bigr] \overline{\nu}(dc),
\nonumber
\\[-8pt]
\\[-8pt]
\nonumber
\overline{\nu}(dc) &=&\frac{1(c>0)}{C_{\nu}}c\nu(dc),\qquad C_{\nu
}=\int
_{\mathbb{R}_+}c\nu(dc).
\end{eqnarray}

We denote
\[
\chi= \biggl( 1-\frac{\alpha}{2} \biggr) \zeta+\frac{\zeta-1}{2}
\quad\mbox{and} \quad\rho=
\frac{\chi}{h}.
\]
Since $1-\frac{\alpha}{2}<h$, there exists $\zeta\in(1,\rho^{-1})$
with $\rho\in(0,1)$. We fix
such a $\zeta$.
We define now
\[
{\gamma} _{t}(x,y):=t^{\chi/ h }g_{t}(x,y)
\]
and we notice that by (\ref{eq:TP})
\[
\bigl\llvert \hat{\theta}_{t}(x,y)\bigr\rrvert \phi_{t}^{y}(x)
\leq Cg_{t}(x,y)=Ct^{-\chi/ h}{\gamma} _{t}(x,y)=Ct^{-\rho}{
\gamma} _{t}(x,y).
\]
We also define $G_t(x,y):=\mathbb{E} [ q_{\sigma
^{2}(y)V_{t}}(x-y) ]=p_{t}^{y}(x,y) $, and we use
Lemma~\ref{lem:estforq} in order to define $\gamma^3$
\begin{eqnarray*}
\bigl\llvert \partial_{x}p_{t}^{y}(x,y)\bigr
\rrvert &\leq&\mathbb{E} \biggl[ \frac
{\llvert
x-y\rrvert  }{\sigma^{2}(y)V_{t}}q_{\sigma^{2}(y)V_{t}}(x-y) \biggr]
\leq C\mathbb{E}%
 \bigl[ V_{t}^{-{1} /2}q_{\sigma^{2}(y)V_{t}}(x-y)
\bigr] \\
&=:&Ct^{-{1}/{(2h)}}\gamma^3_t(x,y).
\end{eqnarray*}
With these definitions, we need to check that \eqref{H2} and \eqref
{h2a} holds.
We verify the former as the latter is similar to \eqref{eq:rost} if one
uses \eqref{eq:estV} at the end of the calculation.
To verify \eqref{H2}, it is enough to prove that for $n\in\mathbb
{N},\delta_{i}>0,i=1,\ldots,n$
%
\begin{equation}
\label{eq:P6} \sup_{y_{0},y_{n+1}\in\mathbb{R}^d}\int dy_{1}\cdots\int
dy_{n}\prod_{i=0}^{n}
\gamma_{\delta_{i}}(y_{i},y_{i+1})^{\zeta}\leq
\frac{C^{n}}{(\delta_{1}+\cdots+\delta_{n})^{{1}/{(2h)}}}, 
\end{equation}
where $C$ is a constant which depends on $\zeta$, $p$, $h$ and $s_{\ast
}$ which
appear in Hypothesis~\ref{(P2)} and in (\ref{P2})$.$ Notice first that
for every $a>0$
and $x\in\mathbb{R}$ one has for a positive constant $C$,
\[
\bigl( q_{a}(x) \bigr) ^{\zeta}=Ca^{-{(\zeta-1)}/{2}}q_{{a}/\zeta}(x).
\]
Using H\"{o}lder's
inequality and the definition of $\chi$, we obtain
\begin{eqnarray*}
g_{\delta_{i}}(y_{i},y_{i+1})^{\zeta} &\leq&\int
_{\mathbb{R}_+}\mathbb {E} \bigl[ {%
V_{\delta_{i}}^{-(1-\alpha/2)\zeta}}
q_{({\overline{\sigma
}^{2}}/{2})%
(c+V_{\delta_{i}})}(y_{i}-y_{i+1})^{\zeta} \bigr]
\overline{\nu}(dc)\\
& \leq& C\int_{\mathbb{R}_+}\mathbb{E} \bigl[
{V_{\delta_{i}}^{-\chi}} q_{({\overline{\sigma}^{2}}/{2(\zeta)})(c+V_{\delta
_{i}})}(y_{i}-y_{i+1})%
 \bigr] \overline{\nu}(dc).
\end{eqnarray*}
We consider $(V_{t}^{i})_{t\geq0},i=1,\ldots,n$ to be independent
copies of $(V_{t})_{t\geq0}$ and we write
\begin{eqnarray*}
&&\int dy_{1}\cdots\int dy_{n}\prod
_{i=0}^{n}g_{\delta
_{i}}(y_{i},y_{i+1})^{\zeta}
\\
&&\qquad\leq C^{-n}\mathbb{E} \Biggl[ \prod_{i=1}^{n}
\bigl(V_{\delta
_{i}}^{i}\bigr)^{-\chi}\int\overline{
\nu}(dc_{1})\cdots\int\overline{\nu }(dc_{n})\int
dy_{1}\cdots\int dy_{n}\\
&&\hspace*{115pt}\qquad{}\times \prod_{i=1}^{n}q_{({\overline{\sigma
}^{2}}/{(2\zeta)})%
(c_{i}+V_{\delta_{i}}^{i})}(y_{i}-y_{i+1})
\Biggr]
\\
&&\qquad=C^{-n}\mathbb{E} \Biggl[ \prod_{i=1}^{n}
\bigl(V_{\delta_{i}}^{i}\bigr)^{-\chi
}\int \overline{
\nu}(dc_{1})\cdots\\
&&\hspace*{43pt}\qquad{}\times\int\overline{\nu}(dc_{n})q_{
({\overline{%
\sigma}^{2}}/{(2\zeta)})\sum_{i=1}^{n}(c_{i}+V_{\delta
_{i}}^{i})}(y_{0}-y_{n+1})
\Biggr]
\\
&&\qquad\leq C^{-n} \Biggl( \mathbb{E} \Biggl[ \prod
_{i=1}^{n}\bigl(V_{\delta
_{i}}^{i}
\bigr)^{-2\chi\zeta} \Biggr] \Biggr) ^{1/2} \\
&&\qquad\quad{}\times\biggl( \mathbb{E} \biggl[
\int\overline{\nu}(dc_{1})\cdots\int\overline{\nu}%
(dc_{n})q_{({\overline{\sigma}^{2}}/{(2\zeta)})\sum_{i=1}^{n}(c_{i}+V_{
\delta_{i}}^{i})}^{2}(y_{0}-y_{n+1})
\biggr] \biggr) ^{1/2}.
\end{eqnarray*}
Notice that $V$ is a L\'{e}vy processes, therefore, $\sum_{i=1}^{n}V_{\delta
_{i}}^{i}$ has the same law as $V_{\delta_{1}+\cdots+\delta_{n}}$ so
\begin{eqnarray} \label{eq:estV}
\nonumber
&&\mathbb{E} \biggl[ \int\overline{\nu}(dc_{1})\cdots\int
\overline{\nu }(dc_{n})q_{%
({\overline{\sigma}^{2}}/{(2\zeta)})\sum_{i=1}^{n}(c_{i}+V_{\delta
_{i}}^{i})}^{2}(y_{0}-y_{n+1})
\biggr]
\\
&&\qquad=\mathbb{E} \biggl[ \int\overline{\nu}(dc_{1})\cdots
 \int\overline{
\nu}(dc_{n})q_{%
({\overline{\sigma}^{2}}/{(2\zeta)})(\sum_{i=1}^{n}c_{i}+V_{\delta
_{1}+\cdots+\delta_{n}})}^{2}(y_{0}-y_{n+1})
\biggr]
\nonumber
\\[-8pt]
\\[-8pt]
\nonumber
&&\qquad\leq\mathbb{E} \biggl[ \int \overline{\nu}(dc_{1})\cdots\int
\overline{\nu}(dc_{n})%
\frac{2\zeta}{\overline{\sigma}^{2}(\sum_{i=1}^{n}c_{i}+V_{\delta
_{1}+\cdots+\delta_{n}})} \biggr]
\\
& &\qquad\leq\frac{2\zeta}{\overline{\sigma
}^{2}}\mathbb{E} \biggl[ \frac{1}{V_{\delta
_{1}+\cdots+\delta_{n}}} \biggr] \leq
\frac{C\zeta}{\overline{\sigma
}^{2}(\delta_{1}+\cdots+\delta_{n})^{{1}/h}}.\nonumber
\end{eqnarray}
The last inequality is a consequence of (\ref{P2}). Again by (\ref{P2}),
\[
\Biggl( \mathbb{E} \Biggl[ \prod_{i=1}^{n}
\bigl(V_{\delta_{i}}^{i}\bigr)^{-2\chi
\zeta} \Biggr] \Biggr)
^{{1}/2}=\prod_{i=1}^{n} \bigl(
\mathbb{E} \bigl[ \bigl(V_{\delta
_{i}}^{i}\bigr)^{-2\chi\zeta}
\bigr] \bigr) ^{{1}/2}\leq C^{n}\prod
_{i=1}^{n}\delta_{i}^{-{\chi\zeta}/h}
\]
so (\ref{eq:P6}) is proved.
\end{pf}

We give now a probabilistic representation for the density of the
solution of \eqref{enu}. We consider the Poisson process $%
J$ of parameter $\lambda=1$ with jump times $\{\tau_i;i\in\mathbb{N}\}$,
and a sequence of i.i.d. standard normal random variables $(\Delta
_{j})_{j\in\mathbb{N}}$.
First, note that using the mean value theorem, we can
rewrite (\ref{eq:thetajumps}) as
\begin{eqnarray*}
&&\hat{\theta}_{t}(x,y)\phi_{t}^{y}(x)=C_{\nu}
\int_{\mathbb{R}_+}\int_{\underline{%
\sigma}^{2}}^{\overline{\sigma}^{2}}1_{\{\sigma^{2}(x)\wedge\sigma
^{2}(y)\leq
u\leq\sigma^{2}(x)\vee\sigma^{2}(y)\}}\operatorname{sgn}_{\sigma}(x,y)
\\
&&\hspace*{125pt}{}\times \mathbb {E} \bigl[ \partial _{t}q_{uc+\sigma^{2}(y)V_{t}}(x-y) \bigr] \,du\,
\overline{\nu}(dc).
\end{eqnarray*}
Here, we define
\[
\operatorname{sgn}_{\sigma}(x,y)= \cases{ %
1, &\quad $\mbox{if }
\sigma^{2}(x)>\sigma^{2}(y)$,
\vspace*{2pt}\cr
-1, &\quad $ \mbox{if }\sigma^{2}(x)\leq\sigma^{2}(y).$}
\]
Therefore, we have that if we consider $U_i\sim \operatorname{Unif}[\underline{\sigma
}^{2},%
\overline{\sigma}^{2}]$ an i.i.d. sequence of random variables independent
of all other random variables, we can
represent $\hat{\theta}_{t}(x,y)\phi_{t}^{y}(x)$ as
\begin{eqnarray*}
&&\hat{\theta}_{t}(x,y)\phi_{t}^{y}(x)\\
&&\qquad=
\frac{C_{\nu}}{2} \bigl( \bar {\sigma}%
^{2}-\underline{
\sigma}^{2} \bigr)\\
&&\qquad\quad{}\times \int_{\mathbb{R}}\operatorname{sgn}_{\sigma
}(x,y)
\mathbb{E} \bigl[ 1_{\{\sigma^{2}(x)\wedge\sigma^{2}(y)\leq U\leq\sigma^{2}(x)\vee
\sigma
^{2}(y)\}}\\
&&\qquad\quad\hspace*{81pt}{}\times h_{Uc+\sigma^{2}(y)V_{t}}^{1,1}(x,y)q_{Uc+\sigma
^{2}(y)V_{t}}(x-y)
\bigr] \overline{\nu}(dc).
\end{eqnarray*}
Here, $h^{1,1}$ is the Hermite polynomial defined in (\ref{D12}). In this
case, the approximating Markov chain is defined as $Y_{0}^{\ast,\pi}(y)=y$
\[
Y_{\tau_{i+1}}^{\ast,\pi}(y)=Y_{\tau_{i}}^{\ast,\pi}(y)+\Delta
_{i} \bigl( U_{i}Z_{i}+\sigma^{2}
\bigl(Y_{\tau_{i}}^{\ast,\pi}(y)\bigr) (V_{\tau
_{i+1}}-V_{\tau
_{i}})
\bigr) ^{{1}/2}.
\]
The corresponding weight is given by
\begin{eqnarray*}
&& \Gamma_{T}^{\pi}(y)= \biggl( \frac{C_{\nu}}{2} \bigl(
\bar {\sigma}^{2}-%
\underline{\sigma}^{2} \bigr)
\biggr) ^{J_{T}}\\
&&\hspace*{42pt}{}\times \prod_{i=1}^{J_{T}}\operatorname{sgn}_{\sigma}
\bigl(Y_{\tau_{i+1}}^{\ast,\pi
}(y),Y_{\tau
_{i}}^{\ast,\pi}(y)
\bigr)\\
&&\hspace*{68pt}{}\times 1_{\{\sigma^{2}(Y_{\tau_{i+1}}^{\ast,\pi
}(y))\wedge\sigma
^{2}(Y_{\tau_{i}}^{\ast,\pi}(y))\leq U_{i}\leq\sigma^{2}(Y_{\tau
_{i+1}}^{\ast,\pi
}(y))\vee\sigma^{2}(Y_{\tau_{i}}^{\ast,\pi}(y))\}}
\\
&&\hspace*{68pt}{}\times h_{U_{i}Z_{i}+\sigma
^{2}(Y_{\tau_{i}}^{\ast,\pi}(y))V_{t}}^{1,1} \bigl( Y_{\tau
_{i+1}}^{\ast,\pi}(y),Y_{\tau_{i}}^{\ast,\pi}(y)
\bigr).
\end{eqnarray*}

\begin{corollary}
\label{cor:75}
Under the hypothesis of Theorem~\ref{th:74}, we have
\begin{eqnarray*}
p_{T}(x,y)&=&e^{T}\mathbb{E} \bigl[ p_{T-\tau_{T}}^{Y_{\tau_{T}}^{\ast
,\pi
}(y)}
\bigl(x,Y_{\tau_{T}}^{\ast,\pi}(y)\bigr)\Gamma_{T}^{\pi}(y)
\bigr] , \label{P7}
\\
\partial_{x}p_{T}(x,y)&=&e^{T}\mathbb{E} \bigl[
\partial_{x}p_{T-\tau
_{T}}^{Y_{\tau
_{T}}^{\ast,\pi}(y)}\bigl(x,Y_{\tau_{T}}^{\ast,\pi}(y)
\bigr)\Gamma_{T}^{\pi
}(y) \bigr]. \label{P8}
\end{eqnarray*}
\end{corollary}
%
\subsubsection{Examples of L\'{e}vy measures}

We conclude this section with two examples of L\'evy measures that
satisfy Hypotheses \ref{(P1)} and \ref{(P2)}.

\begin{example}
Let $c_{k}=k^{-\rho}$ for some $\rho>1$ and define the discrete
measure $\nu(dc)
=\sum_{k=1}^{\infty}\delta_{c_{k}}(dc)$.

We verify that all the hypotheses required in Section~\ref{sec:jumps}
are satisfied in this example. First of all,
we consider Hypothesis \ref{(P1)}. Clearly, $\nu(\mathbb{R}_+)=\infty
$, and if $\rho>1$ then
$\int c\,d\nu(c)<\infty$.

Now we verify Hypothesis \ref{(P2)}. One has $\eta_{\nu
}(a)=\operatorname{card}\{k\dvtx c_{k}>a\}=[a^{-{1}/\rho}]-1(a^{-{1}/\rho}\in
\mathbb{N})$ for $a>0$. We
define $\eta_{\nu}^{\prime}(a)=a^{-{1}/\rho}$. Then clearly
$\eta_{\nu
}^{\prime}$ satisfies the Hypothesis \ref{(P2)} with $h=\frac{1}\rho$.
Furthermore, $s^{-{1}/\rho}\eta_{\nu}(\frac{u}s)-\eta_{\nu
}^{\prime}(u)$
converges uniformly to zero as $s\rightarrow\infty$. Then as $\eta
_{\nu
}\leq\eta_{\nu}^{\prime}$ then Hypothesis \ref{(P2)} is verified
for $\eta_{\nu
} $ with $h=\frac{1}\rho$. So we may use Corollary~\ref{cor:75} or
Theorem~\ref{th:74} for equations with $\alpha$-H\"older coefficient
$\sigma$ with
$\alpha>\frac{2(\rho-1)}\rho$ and the
L\'evy measure $\mu_\nu(du)=q_{\nu}(u)\,du$ with
\[
q_{\nu}(u)=\frac{1}{\sqrt{2\pi}}\sum_{k=1}^{\infty}k^{\rho/2}e^{-{k^{\rho
}u^{2}}/2}.
\]
\end{example}

\begin{example} We consider the measure $\nu
(dc)=1_{[0,1]}(c)c^{-(1+\beta)
}\,dc$ with $\frac{1}2 <\beta<1$. Then $\nu(\mathbb{R}_+)=\infty$ and
$\int c\,d\nu
(c)<\infty$. One has $\eta_{\nu}(a)=\frac{1}{\beta}(a^{-\beta}-1)$
for $a\in(0,1)$ so that Hypothesis \ref{(P2)} holds with $h=\beta\in
{}[0,1)$. Therefore, Corollary~\ref{cor:75} or Theorem~\ref{th:74}
can be applied for an $\alpha$-H\"older coefficient $\sigma$ with
$\alpha\in(2(1-\beta), 1)$. One
may also compute
\[
q_{\nu}(u)=\frac{2^{2\beta}}{u^{1-\beta}\sqrt{2\pi}} \int_{
{u^{2}}/{2}}^{\infty}y^{-1-2\beta}e^{-y}\,dy
\]
so we have the following asymptotic behavior around $0$ for $q_\nu$:
\[
\lim_{u\rightarrow0}u^{2\beta}q_{\nu}(u)=
\frac{2^{2\beta}}{%
\sqrt{2\pi}}\int_{0}^{\infty}y^{\beta-1}e^{-y}\,dy<
\infty.
\]
Therefore, the L\'evy measure generated by this example is of
stable-like behavior around $0$.
\end{example}
%
\section{Some conclusions and final remarks}

The parametrix method has been a successful method in the mathematical
analysis of fundamental solutions of PDEs and we wanted to show the
reader the possibility of other directions of possible generalization.
One of them is to use the current set-up to introduce
stochastic processes representing a variety of different operators
which are
generated by a parametrized operator $L^{z}$. Therefore, allowing the
stochastic representation for various nontrivial operators.



The adjoint method we introduced here seems to allow for the analysis
of the regularity of the density requiring H\"older continuity of the
coefficients through an explicit expression of the density.

Finally, the stochastic representation can be used for simulation
purposes. In that case, the variance of the estimators explode due to
the instability of the weight function $\theta$ in the forward method
or $\hat{\theta}$ in the backward method.
In fact, the representations presented here have a theoretical
infinite variance although the mean is finite.
In that respect, the way that the Poisson process and the exponential
jump times appear maybe considered somewhat arbitrary. In fact, one can
think of various other representations which may lead to variance
reduction methods. Preliminary simulations show that different
interpretations of the time integrals in the parametrix method may lead
to finite variance simulation methods.
Many of these issues will be taken up in future work.

\begin{appendix}
\section*{Appendix}\label{app}

\subsection{On some Beta type coefficients}

\label{sec:8.1}
For $t_0\in\mathbb{R}$, $a\in[0,1)$, $b>-1$ and $n\in\mathbb{N}$, define
\[
c_{n}(t_{0},a,b):=\int_{0}^{t_{0}}\,dt_{1}
\cdots%
\int_{0}^{t_{n-1}}\,dt_{n}t_{n}^{b}
\prod_{j=0}^{n-1}(t_{j}-t_{j+1})^{-a}.
\]

\begin{lemma}
\label{lem:app1}Let $a\in[0,1)$ and $b>-1$. Then we have
\[
c_{n}(t_{0},a,b)\leq t_{0}^{b+n(1-a)}
\frac{\Gamma(1+b)\Gamma
^{n}(1-a)}{%
[1+b+n(1-a)]!} \qquad \mbox{for } n\geq\frac{1-b}{1-a}.
\]
In particular, for $b=0$,
%
\renewcommand{\theequation}{\Alph{section}.\arabic{equation}}
\begin{equation}
\qquad c_{n}(t_{0},a):=c_{n}(t_{0},a,b)\leq
t_{0}^{n(1-a)}\frac{\Gamma
^{n}(1-a)}{%
[1+n(1-a)]!} \qquad\mbox{for } n\geq(1-a)^{-1}.
\label{A6bis}
\end{equation}
\end{lemma}

\begin{pf}Let $b>-1$ and $0\le a<1$ and use the change of variable $%
s=ut$ so that
\[
\int_{0}^{t}(t-s)^{-a}s^{b}\,ds=t^{b+1-a}%
\int_{0}^{1}(1-u)^{-a}u^{b}\,du=t^{b+1-a}B(1+b,1-a),
\]
where $B(x,y)=\int_{0}^{1}t^{x-1}(1-t)^{y-1}\,dt$ is the standard Beta
function and $b+1-a>-1$. Using this repeatedly, we obtain
\begin{eqnarray*}
c_{n}(t_{0},a,b)&=&t_{0}^{b+n(1-a)}%
\prod_{i=0}^{n-1}B\bigl(1+b+i(1-a),1-a
\bigr)\\
&=&t_{0}^{b+n(1-a)}\frac{\Gamma
(1+b)\Gamma
^{n}(1-a)}{\Gamma(1+b+n(1-a))}.
\end{eqnarray*}
The last equality being a consequence of the identity $B(x,y)=\frac
{\Gamma
(x)\Gamma(y)}{\Gamma(x+y)}$. The function $\Gamma(x)$ is increasing
for $%
x\geq2$ so the result follows. Letting $b=0$, we get \eqref{A6bis}.
\end{pf}
%
\subsection{Some properties of Gaussian type kernels}
\label{sec:Gaussian}

In this section, we introduce some preliminary estimates concerning Gaussian
kernels. We consider a $d$ dimensional square symmetric nonnegative
definite matrix $a$. We assume that $0<\underline{a}I\leq a\leq
\overline{a}%
I $ for $\underline{a}$, $\overline{a}\in\mathbb{R}$ and we define
$\rho_{a}:=%
\frac{\overline{a}}{\underline{a}}$. The Gaussian density of mean zero
and covariance matrix $%
a$ is denoted by
\[
q_{a}(y)=\frac{1}{(2\pi)^{{d}/2}\sqrt{\det a}}\exp \biggl(-\frac
{1}{2} \bigl
\langle a^{-1}y,y \bigr\rangle \biggr).
\]
For a strictly positive real number $\lambda$, we abuse the notation,
denoting by $q_{\lambda}(y)\equiv q_{a}(y)$ for $a=\lambda I$ where
$I$ is
the identity matrix. In particular, $q_{1}$ is the standard Gaussian kernel
on $\mathbb{R}^d$. Then we have the following immediate inequalities:
%
\renewcommand{\theequation}{\Alph{section}.\arabic{equation}}
\begin{eqnarray}\label{G1}
\mathrm{(i)}&&\quad q_{s}(y) \leq \biggl(\frac{t}{s} \biggr)^{{d}/2}q_{t}(y)\qquad
\forall s<t,
\nonumber
\\[-8pt]
\\[-8pt]
\nonumber
\mathrm{(ii)}&& \quad \rho_{a}^{-{d}/2}q_{\underline{a}}(y) \leq
q_{a}(y)\leq\rho _{a}^{{d}/2}q_{\overline{a}}(y).
\nonumber
\end{eqnarray}
We define for $a\in\mathbb{R}^{d\times d}$, the Hermite polynomials in
$\mathbb{R}^d$ as
\[
H_{a}^{i}(y)=-\bigl(a^{-1}y
\bigr)^{i},\qquad H_{a}^{i,j}(y)=\bigl(a^{-1}y
\bigr)^{i}\bigl(a^{-1}y\bigr)^{j}-
\bigl(a^{-1}\bigr)^{i,j}. \label{G2}
\]
Direct computations give
%
\begin{equation}
\partial_{i}q_{a}(y)=H_{a}^{i}(y)q_{a}(y),\qquad
\partial _{i,j}^{2}q_{a}(y)=H_{a}^{i,j}(y)q_{a}(y).
\label{G3}
\end{equation}
We will use the following basic estimates.

\begin{lemma}
\label{lem:estforq}For $\alpha\in{}[0,1]$, we have for all
$i,j\in
\{1,\ldots,d\}$, $y\in\mathbb{R}^d$ and $t>0$,
%
\begin{eqnarray}\label{G4}
\mathrm{(i)}&&\quad \llvert y\rrvert ^{\alpha}\bigl\llvert \partial _{i,j}^{2}q_{ta}(y)
\bigr\rrvert\leq C_{a}\frac{1}{t^{1-\alpha/
2}}q_{{t%
\overline{a}}/2}(y) \quad\mbox{and}
\nonumber
\\[-8pt]
\\[-8pt]
\nonumber
\mathrm{(ii)}&&\quad \llvert y\rrvert ^{\alpha}\bigl\llvert \partial_{i}q_{ta}(y)
\bigr\rrvert \leq C_{a}^{\prime}\frac{1}{t^{{(1-\alpha)}/{2}}}q_{{t\overline{a}}/2}(y)
\nonumber
\end{eqnarray}
with
\[
C_{a}= ( 2\rho_{a} ) ^{{d/}2}{
\underline{a}^{-1}} (4%
\overline{a})^{\alpha/2}(4
\rho_{a}+1),\qquad C_{a}^{\prime}={%
\underline{a}^{-1}} ( 4\overline{a} ) ^{{(1+\alpha)}/2} ( 2\rho
_{a} ) ^{{d}/2}.
\]
\end{lemma}

\begin{pf}We have
\[
\bigl\llvert H_{ta}^{i,j}(y)\bigr\rrvert \leq
\frac{\llvert  y\rrvert
^{2}}{\underline{a}%
^{2}t^{2}}+\frac{1}{\underline{a}t}
\]
so that
\[
\llvert y\rrvert ^{\alpha}\bigl\llvert \partial_{i,j}^{2}q_{ta}(y)
\bigr\rrvert \leq \frac{1}{\underline{a}t^{1-\alpha/2}} \frac{\llvert  y\rrvert
^{\alpha}%
}{t^{\alpha/2}} \biggl( \frac{\llvert  y\rrvert  ^{2}}{\underline{a}t}
+1 \biggr) q_{ta}(y).
\]
We use (\ref{G1}) and we obtain
\[
q_{ta}(y)\leq\rho_{a}^{{d}/2}q_{t\overline{a}}(y)=
( 2\rho _{a} ) ^{{d}/2}\exp \biggl(-\frac{\llvert  y\rrvert  ^{2}}{4t\overline
{a}}
\biggr)q_{{t
\overline{a}}/2}(y).
\]
We may find a constant $c_{\alpha}$ such that $v^{\lambda}e^{-v}\leq
c_{\alpha}$ for every $0\leq\lambda\leq2+\alpha$. Using this inequality
twice with $\lambda=%
\frac{2+\alpha}{2}$ and $\lambda=\frac{\alpha}{2}$ and for $v=\frac
{1}{4t\overline{a}}\llvert  y\rrvert  ^{2}$, we obtain
\[
\llvert y\rrvert ^{\alpha}\bigl\llvert \partial_{i,j}^{2}q_{ta}(y)
\bigr\rrvert \leq \frac{ ( 2\rho_{a} ) ^{{d}/2}}
{\underline{a}t^{1-\alpha/2}} (4%
\overline{a})^{\alpha/2}(4
\rho_{a}+1)q_{{t\overline{a}}/2}(y).
\]
The proof of (ii) is similar.
\end{pf}

\begin{lemma}
\label{lem:ineqgauss}Let $0<s'\le s''$ and $y\in\mathbb{R}$ then
\[
\bigl\llvert q_{s^{\prime\prime}}(y)-q_{s^{\prime}}(y)\bigr\rrvert \le \sqrt{
\frac{2s^{\prime\prime}}{s^{\prime}}}q_{s^{\prime
\prime}}(y) \biggl( \frac{s^{\prime\prime}-s^{\prime}}{s^{\prime
}} \biggr) .
\]
\end{lemma}
\begin{pf} Using the fact that $q$ solves the heat equation, we have
using \eqref{G1}
and \eqref{G4}
%
\begin{eqnarray}\label{eq:A}
\nonumber
\bigl\llvert q_{s^{\prime\prime}}(y)-q_{s^{\prime}}(y)\bigr\rrvert &
\leq &\int_{s^{\prime}}^{s^{\prime\prime}}\bigl\llvert \partial
_{s}q_{s}(y)\bigr\rrvert \,ds=\frac{1}{2}\int
_{s^{\prime}}^{s^{\prime\prime}}\bigl\llvert \partial
_{a}^{2}q_{s}(y)\bigr\rrvert \,ds\leq C\int
_{s^{\prime}}^{s^{\prime\prime
}}\frac{1}{s}q_{{s}/2}(y)\,ds
\\
&\leq&\sqrt{%
\frac{s^{\prime\prime}}{s^{\prime}}}q_{{s^{\prime\prime
}}/2}(y)\int
_{s^{\prime}}^{s^{\prime\prime}}\frac{1}{s}\,ds\leq\sqrt {
\frac{%
2s^{\prime\prime}}{s^{\prime}}}q_{s^{\prime\prime}}(y)\ln \biggl( \frac{%
s^{\prime\prime}}{s^{\prime}} \biggr)
\\
\nonumber
&=&\sqrt{\frac{2s^{\prime\prime}}{s^{\prime}}}q_{s^{\prime
\prime
}}(y)\ln \biggl( 1+ \biggl(
\frac{s^{\prime\prime}}{s^{\prime}}-1 \biggr) \biggr) \leq\sqrt{\frac{2s^{\prime\prime}}{s^{\prime}}}q_{s^{\prime
\prime}}(y)
\biggl( \frac{s^{\prime\prime}-s^{\prime}}{s^{\prime
}} \biggr) .
\end{eqnarray}
From here, the result follows.
\end{pf}
\end{appendix}

\section*{Acknowledgments}
A. Kohatsu-Higa would like to thank L. Li for comments on an earlier
version of this article.





\printaddresses
\end{document}